\documentclass[final]{amsart}
\usepackage{amssymb, amsthm, mathrsfs}
\usepackage{biblatex}
\usepackage{enumitem}
\usepackage{float}
\usepackage{hyperref}
\usepackage{tikz, tikz-cd}

\theoremstyle{plain}
\newtheorem{conjecture*}{Conjecture}
\newtheorem{proposition*}[conjecture*]{Proposition}
\newtheorem{theorem*}[conjecture*]{Theorem}

\newtheorem{corollary}{Corollary}[section]
\newtheorem{lemma}[corollary]{Lemma}
\newtheorem{proposition}[corollary]{Proposition}
\newtheorem{theorem}[corollary]{Theorem}

\newtheorem{conjecture}[corollary]{Conjecture}

\theoremstyle{definition}

\newtheorem{example}[corollary]{Example}
\newtheorem{remark}[corollary]{Remark}

\newcommand{\qbinom}{\genfrac{[}{]}{0pt}{}}

\newcommand{\calB}{\mathcal B}

\newcommand{\calH}{\mathcal H}

\newcommand{\calP}{\mathcal P}

\newcommand{\calR}{\mathcal R}

\newcommand{\bbN}{\mathbb N}
\newcommand{\bbQ}{\mathbb Q}
\newcommand{\bbR}{\mathbb R}
\newcommand{\bbZ}{\mathbb Z}

\newcommand{\frakg}{\mathfrak g}
\newcommand{\frakh}{\mathfrak h}

\newcommand{\fraksl}{\mathfrak{sl}}

\newcommand{\ula}{{\underline a}}

\newcommand{\uli}{{\underline i}}

\newcommand{\te}{\tilde e}
\newcommand{\tf}{\tilde f}

\newcommand{\cf}{\check{f}}
\newcommand{\ce}{\check{e}}
\newcommand{\cepsilon}{\check{\epsilon}}
\newcommand{\cphi}{\check{\varphi}}

\DeclareMathOperator{\id}{id}

\DeclareMathOperator{\supp}{supp}
\DeclareMathOperator{\trace}{tr}
\DeclareMathOperator{\weight}{wt}

\DeclareMathOperator{\ad}{ad}

\DeclareMathOperator{\Hom}{Hom}

\DeclareMathOperator{\pr}{pr}

\title{Base change between canonical basis and path basis}
\author{Felix Röhrich}
\address{Lehrstuhl für Algebra und Darstellungstheorie, RWTH Aachen University, Pontdriesch 10-16, 52062 Aachen, Germany}
\email{roehrich.math@pm.me}

\begin{document}

\begin{abstract}
    Almost 30 years ago Littelmann conjectured that there is a unitriangular base change between the dual canonical basis and the path vector basis for all symmetrizable Kac-Moody Lie algebras.
    In this paper we give a positive answer to the conjecture for the case of classical Lie algebras and type $G_2$, when the path model is given by standard Young tableaux.
\end{abstract}

\maketitle

\section{Introduction}

Almost 30 years ago, Littelmann \cite{Littelmann:ContractingModules-StandardMonomialTheory} conjectured that there is a unitriangular base change between the dual canonical basis and the path vector basis for all symmetrizable Kac-Moody Lie algebras.
However, little progress has been made since then.
This is primarily for two reasons:
First, the path vector basis is incompatible with the crystal structure of the canonical basis.
This can be seen from the fact that the path basis is not perfect in the sense of Baumann, Kamnitzer and Knutson \cite{BaumannKamnitzerKnutson:MirkovicVilonenBasis-DuistermaatHeckmanMeasures} (the definition of \emph{perfect} by Berenstein and Kazhdan \cite{BerensteinKazhdan:Geometric-UnipotentCrystals-II} is identical up to normalization).
Second, the canonical basis is incompatible with the Frobenius splitting in general \cite{Baumann:CanonicalBasis-QuantumFrobeniusMorphism,Roehrich:DualImaginaryVectors-TightMonomialCones-QuantumFrobeniusMorphism}.
However, the definition of the path vectors makes use of the Frobenius splitting.
Thus, to prove the conjecture a new perspective is required.

Let $\frakg$ be symmetrizable Kac-Moody Lie algebra and $U_q(\frakg)$ its quantized enveloping algebra.
Consider a dominant weight $\lambda$ and its associated integrable highest weight module $V(\lambda)$.
Recall that the combinatorics of the crystal base $B(\lambda)$ can be realized through Littelmann's path model \cite{Kashiwara:Similarity-CrystalBases}.
For our purposes we will only consider Lakshmibai-Seshadri paths and their concatenations.
Kashiwara's operators $\te_i, \tf_i$ on $B(\lambda)$ translate into the so-called root operators which are denoted by same symbols.
For a path $\pi\colon [0,1] \to P \otimes_\bbZ \bbR$, where $P$ is the weight lattice, the action of root operators is determined by the global minima of the function $\langle h_i, \pi(t)\rangle$, where $h_i$ is the simple coroot associated to the simple root $\alpha_i$.
The combinatorial structure of the path vectors, however, does not depend on some global minimum but on the first local minimum.
Now, it is easy to modify the definition of the root operators to obtain new operators $\ce_i, \cf_i$ which we will refer to as \emph{modified root operators}.
We can also define functions $\cepsilon_i, \cphi_i, \weight$ to obtain a structure that is very close to a crystal in the sense of Kashiwara \cite{Kashiwara:CrystalBase-Littelmann-DemazureCharacterFormula}.
The main obstruction is that the new structure is not compatible with the tensor product rule of crystals.

Recall that any realization of $B(\lambda)$ in terms of the path model is uniquely determined by a \emph{dominant path} $\eta$, a path contained entirely in the dominant Weyl chamber, ending in $\lambda$.
Moreover, if $\eta$ and $\eta'$ are dominant paths of the same weight $\lambda$, then $B(\eta) \cong B(\eta')$ as crystals \cite{Littelmann:Paths-RootOperators-RepresentationTheory}.
Hence, the action of the root operators $\te_i, \tf_i$ is independent of the choice of the dominant path $\eta$.
However, this is not true for the modified root operators $\ce_i$ and $\cf_i$:
If $\pi \in B(\eta)$ corresponds to the path $\pi' \in B(\eta')$ under the crystal isomorphism $B(\eta) \cong B(\eta')$, then $\ce_i\pi, \cf_i\pi$ need not correspond to $\ce_i\pi',\cf_i\pi'$.
In other words, the definition of the modified root operators depends on the choice of the dominant path $\eta$.
Since we work with a fixed path $\eta$, we will omit the dependency of the operators $\ce_i,\cf_i$ on $\eta$ in the notation.
Still, the reader should keep in mind that we are discussing not a single new combinatorial structure on the canonical basis but multiple at once.

Our observations show that the action of the operator $\cf_i$ is a combinatorial shadow of left action of the Chevalley generator $f_i \in U_q^-(\frakg)$ on an element of the canonical basis of $V(\lambda)$, which is expressed in the rather technical Conjecture \ref{conjecture:f-action}.
Under mild assumptions this conjecture can be simplified to the following formula:
\[
        f_i^{(n)} G_\lambda(b_\pi) = \qbinom{\cepsilon_i(\pi)+n}{n}_i G_\lambda(b_{\cf_i^n\pi}) + \sum_{\pi'\triangleleft \cf_i^n\pi} F_{\pi,\pi'}^{i,n} G_\lambda(b_{\pi'}).
\]
Here $b_\pi \in B(\lambda)$ is the element of the crystal base corresponding to the path $\pi \in B(\eta)$; the partial order $\trianglerighteq$ is introduced in Section \ref{section:partial-orders}.
This formula on the matrix coefficients has a striking resemblance to the one by Kashiwara (see (\ref{equation:f-left-multiplication-canonical-basis})), which is used in many proofs involving the canonical basis.
At this point there is no apparent connection between the two formulas.
We have the following consequence.

\begin{proposition*}[Proposition \ref{proposition:base-change-path-monomial-basis-canonical-basis}]
    Assume that Conjecture \ref{conjecture:f-action} holds for the dominant path $\eta$ ending in the dominant weight $\lambda$.
    Then, there is a unitriangular base change between the path monomial basis $\calP(\eta)$ and the canonical basis of $V(\lambda)$ for any linearisation of the partial order $\trianglerighteq$ on $B(\eta)$.
\end{proposition*}

The \emph{path monomial basis} is a monomial basis of $V(\lambda)$ that was introduced by Littelmann \cite{Littelmann:ContractingModules-StandardMonomialTheory} to show that the path vectors form a basis.
In particular, if we prove Conjecture \ref{conjecture:f-action}, then Littelmann's conjecture on the unitriangular base change between path vectors and dual canonical basis follows immediately.
As a special case of the above proposition we were able to prove the following result.

\begin{proposition*}[Proposition \ref{proposition:canonical-basis-classical-fundamental-weight} and Lemma \ref{lemma:conjecture-g2-fundamental-weights}]
    \label{proposition*:canonical-basis-fundamental-weight}
    Let $\lambda$ be a fundamental weight of classical type or a fundamental weight in $G_2$.
    Then, the canonical basis of $V(\lambda)$ coincides with the path monomial basis.
\end{proposition*}

The result for fundamental weights of classical type was inspired by Marsh's results \cite{Marsh:CanonicalBasis-StandardMonomials} on the standard monomial basis.
In fact, we will see in Section \ref{section:canonical-basis-fundamental-weights} that definition of path monomial basis and the standard monomial basis coincide for fundamental weights of classical type.

For a general dominant weight $\lambda$ consider the decomposition $\lambda = m_1\varpi_1 + \dots + m_N\varpi_N$.
Following Littelmann \cite{Littelmann:PlacticAlgebra-SemisimpleLieAlgebras} we will refer to the elements of $B(\pi_{m_1\varpi_1} * \dots * \pi_{m_N\varpi_N})$ as \emph{standard Young tableaux}.
Here $*$ denotes the concatenation of paths and $\pi_{m_k\varpi_k}$ denotes the straight-line path $t \mapsto t(m_k\varpi_k)$.
Now, the main result of this paper can be stated as follows.

\begin{theorem*}[Theorem \ref{theorem:base-change-path-monomial-basis-canonical-basis-classical-g2}]
    \label{theorem*:base-change}
    There is a unitriangular base change between canonical basis and path monomial basis for classical Lie algebras and for type $G_2$, when the path model is given by standard Young tableaux.
\end{theorem*}

The paper is structured as follows.
In Section \ref{section:preliminaries} we recall the definition of quantized enveloping algebras, canonical basis and path model.
In particular, we will revisit the quasi-$\calR$-matrix and the canonical basis of tensor products of simple integrable $U_q(\frakg)$-modules.
We recommend that the reader reviews this section even if they are familiar with the topics, since some careful choices are necessary to make all conventions align.

In Section \ref{section:base-change} we present the results of this paper.
We begin by introducing the \emph{modified root operators $\ce_i, \cf_i$} and related functions (Section \ref{section:modified-root-operators}), some important partial orders on the path model (Section \ref{section:partial-orders}) and the path monomial basis (Section \ref{section:bases-associated-to-paths}).
Afterwards, we present the central conjectures of this paper (Section \ref{section:conjectures-on-canonical-basis}).
Then, we compare the canonical basis and the path monomial basis for fundamental weights (Section \ref{section:canonical-basis-fundamental-weights}), in particular we obtain Proposition \ref{proposition*:canonical-basis-fundamental-weight}.
Finally, we prove the conjectures for standard Young tableaux paths for classical Lie algebras and in $G_2$ (Section \ref{section:proof-conjecture-standard-young-tableaux-paths}), and obtain Theorem \ref{theorem*:base-change} as a consequence.

\subsection*{Acknowledgements}
I want to thank Xin Fang for helpful discussions and comments on the manuscript.
I acknowledge funding by the Deutsche Forschungsgemeinschaft (DFG, German Research Foundation) - project number 442047500 - through the Collaborative Research Center ``Sparsity and Singular Structures'' (SFB 1481).

\section{Preliminaries}
\label{section:preliminaries}

In this paper we mainly follow the notation and conventions of Kashiwara \cite{Kashiwara:CrystalBases-QuantizedUniversalEnvelopingAlgebras} and Kimura \cite{Kimura:QuantumUnipotentSubgroup-DualCanonicalBasis}.
In examples we use the conventions of Bourbaki \cite{Bourbaki:LieGroups-LieAlgebras-4-6}, Plates I - IX for the labeling of Dynkin diagrams and simple roots.

The notation for the path model is based on \cite{Littelmann:Paths-RootOperators-RepresentationTheory}.
In this paper a path will always mean a concatenation of LS paths unless stated otherwise.

\subsection{Quantized enveloping algebras}

Let $\bbQ(q)$ the field of rational functions in the indeterminate $q$ over $\bbQ$.
We denote by $A$ the $\bbQ$-subalgebra of $\bbQ(q)$ consisting of functions which are regular at $q = 0$.

\subsubsection*{Root Datum}
\label{section:definition-root-datum}

Throughout this paper we fix a \emph{root datum} consisting of
\begin{enumerate}
	\item A finite dimensional $\bbQ$-vector space $\frakh$,
	\item a finite index set $I$,
	\item a linearly independent set $\{ \alpha_i \mid i \in I \}$ of $\frakh^*$ (the set of \emph{simple roots}),
	\item a subset $\{ h_i \mid i\in I\}$ of $\frakh$ (the set of \emph{simple coroots}),
	\item a $\bbQ$-valued symmetric bilinear form $(\cdot, \cdot)$ on $\frakh^*$, and
	\item a lattice $P \subseteq \frakh^*$ (the \emph{weight lattice}) together with the natural pairing $\langle \cdot, \cdot\rangle \colon P^\vee \otimes P \to \bbZ$, where $P^\vee := \Hom_\bbZ(P, \bbZ)$ is the dual lattice,
\end{enumerate}
such that for $i\in I$ the following properties are fulfilled:
\begin{enumerate}
	\item $(\langle h_i, \alpha_j\rangle)_{ij}$ is a generalized Cartan matrix,
	\item $(\alpha_i, \alpha_i) \in 2\bbZ_{>0}$,
	\item $\langle h_i, \lambda\rangle = 2(\alpha_i, \lambda)/(\alpha_i, \alpha_i)$ for $\lambda \in P$,
	\item $\alpha_i\in P$ and $h_i\in P^\vee$.
\end{enumerate}

We assume that the root datum is \emph{simply-connected}, that is, there exist $\varpi_i \in P$ such that $\langle h_i, \varpi_j\rangle = \delta_{ij}$ for any $i, j\in I$.

Let $Q = \bigoplus_{i\in I}\bbZ\alpha_i \subseteq P$ be the \emph{root lattice}; we set $Q_+ := \bigoplus_{i\in I} \bbZ_{\ge0} \alpha_i$ and $Q_- := -Q_+$.
For $\xi = \sum_{i\in I} \xi_i\alpha_i \in Q$, define $\trace(\xi) = \sum_{i\in I} \xi_i$.
We call a weight $\lambda \in P$ \emph{dominant} if $\langle h_i, \lambda\rangle \ge 0$ for all $i\in I$, and denote set of dominant weights by $P_+$.
A weight $\lambda$ is \emph{regular} if $\langle h_i, \lambda\rangle \neq 0$ for all $i\in I$.
We define the \emph{Weyl vector} $\rho := \sum_{i\in I} \varpi_i$.

The weight lattice $P$ has a partial order given by $\lambda \succeq \mu$ if $\lambda - \mu \in Q_+$, that is, $\lambda - \mu$ is a $\bbZ_{\ge0}$-linear combination of simple roots.

\subsubsection*{Kac-Moody Lie algebra}
\label{section:definition-kac-moody-lie-algebra}

Let $\frakg$ be the \emph{symmetrizable Kac-Moody Lie algebra} associated with the Cartan matrix $(\langle h_i, \alpha_j\rangle)_{ij}$.
That is to say, $\frakg$ is the $\bbQ$-vector space generated by $e_i$, $f_i$ ($i \in I$) and $\frakh$ with the Lie bracket
\[
    [h, h'] = 0,\quad [h, e_i] = \langle h, \alpha_i\rangle e_i,\quad [h, f_i] = -\langle h, \alpha_i\rangle f_i, \quad [e_i, f_j] = \delta_{ij} h_i,
\]
for $i, j \in I$ and $h, h' \in \frakh$, and the relations for $i \neq j$
\[
    \ad(e_i)^{1-\langle h_i, \alpha_j\rangle} e_j = 0, \quad \ad(f_i)^{1-\langle h_i, \alpha_j\rangle} f_j = 0.
\]

We write $\frakg_i$ for Kac-Moody Lie subalgebra of $\frakg$ generated by $e_i$, $f_i$; it is isomorphic to $\fraksl_2$.

An element $\alpha \in Q$ is called a \emph{root} of $\frakg$ if the \emph{root space}
\[
    \frakg_\alpha := \{ x \in \frakg \mid [h, x] = \langle h, \alpha\rangle x \text{ for all } h \in \frakh \}
\]
is non-trivial.
A root $\alpha$ is called \emph{positive} (\emph{negative}) if $\alpha \in Q_+$ ($\alpha \in Q_-$).

\subsubsection*{Weyl Group}
\label{section:definition-weyl-group}

For $i\in I$ define the \emph{simple reflection} $s_i\colon \frakh^* \to \frakh^*$,
\begin{equation}
    \label{equation:simple-reflection}
    s_i \lambda = \lambda - \langle h_i, \lambda\rangle \alpha_i.
\end{equation}
The group $W$ generated by the $s_i$ ($i\in I$) is called the \emph{Weyl group}.
The Weyl group also acts on $\frakh$ via
\begin{equation}
    s_i h = h - \langle h, \alpha_i\rangle h_i.
\end{equation}
We have $\langle h, s_i \lambda\rangle = \langle s_i h, \lambda\rangle$ for $h\in \frakh$ and $\lambda \in \frakh^*$.

For a finite sequence $\uli = (i_1, \dots, i_r)$ in $I$ we write
\begin{equation}
    s_\uli = s_{i_1} \cdots s_{i_r}.
\end{equation}

Given $w \in W$, its \emph{length} $\ell(w)$ is the minimal integer $r$ such that $w = s_\uli$ for some sequence $\uli = (i_1, \dots, i_r)$ in $I$; such a word $s_\uli = s_{i_1}\cdots s_{i_r}$ is called a \emph{reduced decomposition} of $w$.
The set of all sequences $\uli$ in $I$ such that $s_\uli$ is a reduced decomposition of $w$ is denoted by $\calH(w)$.

For $w, w' \in W$ we write $w > w'$ if $w$ is greater than $w'$ in the Bruhat order.

For a weight $\lambda \in P$, we denote its stabilizer subgroup in $W$ by $W_\lambda$; we have
\[
    W_\lambda = \langle s_i \mid \langle h_i, \lambda\rangle = 0, i \in I \rangle.
\]
Each coset in $W/W_\lambda$ has a unique representative of minimal length; we write $W^\lambda$ for the set of these representatives.

If the root datum is of finite type, then $W$ is a finite group and has unique element $w_0$ of greatest length.

\subsubsection*{Quantized enveloping algebra}
\label{section:definition-quantized-enveloping-algebra}

The \emph{quantized enveloping algebra} $U_q(\frakg)$ is the algebra over the rational function field $\bbQ(q)$ generated by $e_i, f_i$ ($i\in I$) and $q^h$ ($h\in P^\vee$) subject to the following relations for $i, j \in I$ and $h, h' \in P^\vee$:
\begin{enumerate}
    \item $q^0 = 1$, $q^h q^h = q^{h+h'}$,
    \item $q^h e_i q^{-h} = q^{\langle h, \alpha_i\rangle} e_i$, $q^h f_i q^{-h} = q^{-\langle h, \alpha_i\rangle} f_i$
    \item $[e_i, f_j] = \delta_{ij} (t_i - t_i^{-1})/(q_i - q_i^{-1})$, and
    \item $\sum_{k=0}^a (-1)^k e_i^{(k)} e_j e_i^{(a-k)} = \sum_{k=0}^a (-1)^k f_i^{(k)} f_j f_i^{(a-k)} = 0$ for $i \neq j$, $a = 1 - \langle h_i, \alpha_j\rangle$.
\end{enumerate}
Here we set $q_i = q^{(\alpha_i,\alpha_i)/2}$, $t_i = q_i^{h_i}$; $e_i^{(n)} = e_i^n/[n]_i!$ and $f_i^{(n)} = f_i^n/[n]_i!$ are the \emph{divided powers} as usual, $[n]_i = (q_i^n - q_i^{-n})/(q_i - q_i^{-1})$, $[n]_i! := [n]_i [n-1]_i \cdots [1]_i$ and
\[
    \qbinom{n}{k}_i := \frac{[n]_i [n-1]_i\cdots [n-k+1]_i}{[k]_i [k-1]_i \cdots [1]_i} \quad \text{for } n \in \bbZ \text{ and } k \in \bbN.
\]

We will also use the notation
\begin{equation}
    \qbinom{t_i; n}{k} := \prod_{j=1}^k \frac{q_i^{n+1-j} t_i - q_i^{-(n+1-j)} t_i^{-1}}{q_i^j - q_i^{-j}}
\end{equation}
for $n\in \bbZ$ and $k\in \bbZ_{\ge0}$.

For $n\in \bbZ$ we use the convention that $f_i^{(n)} = e_i^{(n)} = 0$ if $n < 0$.

For a sequence $\uli = (i_1, \dots, i_r)$ in $I$ and $\ula = (a_1, \dots, a_r) \in \bbZ^r$, we set
\begin{equation}
    f_\uli^{(\ula)} := f_{i_1}^{(a_1)} \cdots f_{i_r}^{(a_r)} \in U_q^-(\frakg).
\end{equation}

\subsubsection*{Weight grading}
The algebra $U_q(\frakg)$ has a natural $Q$-grading; for $\xi \in Q$ the graded component is
\[
    U_q(\frakg)_\xi = \{ u \in U_q(\frakg) \mid q^h u q^{-h} = q^{\langle h, \xi\rangle} u \text{ for all } h \in P^\vee \}.
\]
If $u \in U_q(\frakg)_\xi$, we say that $\xi$ is the \emph{weight} of $u$ and denote it by $\weight u$.

\subsubsection*{Automorphisms}

The \emph{bar involution} is the $\bbQ$-algebra automorphism $-\colon U_q(\frakg) \to U_q(\frakg)$ given by
\begin{equation}
    \overline{e_i} = e_i,\quad \overline{f_i} = f_i,\quad \overline{q}= q^{-1},\quad \overline{q^h} = q^{-h}.
\end{equation}

The \emph{$*$-involution} is the $\bbQ(q)$-algebra antiautomorphism $*\colon U_q(\frakg) \to U_q(\frakg)$ given by
\begin{equation}
    *e_i = e_i,\quad *f_i = f_i,\quad *q^h = q^{-h}.
\end{equation}

\subsubsection*{Hopf algebra structure}

In this paper we choose the following coproduct:
\begin{equation}
    \label{equation:coproduct}
    \begin{aligned}
        \Delta(q^h) &= q^h \otimes q^h,\\
        \Delta(e_i) &= e_i\otimes t_i^{-1} + 1\otimes e_i, \\
        \Delta(f_i) &= f_i\otimes 1 + t_i\otimes f_i,
    \end{aligned}
\end{equation}

With respect to $\Delta$, the counit $\epsilon\colon U_q(\frakg) \to \bbQ(q)$ is
\begin{equation}
    \epsilon(e_i) = 0,\quad \epsilon(f_i) = 0,\quad \epsilon(q^h) = 1,
\end{equation}
and the antipode $S\colon U_q(\frakg) \to U_q(\frakg)$ is given by
\begin{equation}
    S(e_i) = -e_i t_i,\quad S(f_i) = -t_i^{-1} f_i,\quad S(q^h) = q^{-h}.
\end{equation}

\subsubsection*{Triangular decomposition}
\label{section:triangular-decomposition}

Let $U_q^-(\frakg)$ be the $\bbQ(q)$-subalgebra of $U_q(\frakg)$ generated by the $f_i$ ($i\in I$) and $U_q^+(\frakg)$ the subalgebra generated by the $e_i$ ($i\in I$).
Then, we have the \emph{triangular decomposition} \cite{Lusztig:IntroductionToQuantumGroups}
\[
    U_q(\frakg) \cong U_q^-(\frakg) \otimes_{\bbQ(q)} \bbQ(q)[P^\vee] \otimes_{\bbQ(q)} U_q^+(\frakg).
\]
Here, $\bbQ(q)[P^\vee]$ is the group algebra of $P^\vee$ over $\bbQ(q)$.

\subsubsection*{Integrable modules}
\label{section:integrable-modules}

Let $M$ be a $U_q(\frakg)$-module. For any $\lambda \in P$, we define its \emph{weight space}
\[
    M_\lambda = \{ m \in M \mid q^h m = q^{\langle h, \lambda\rangle} m \text{ for all } h \in P^\vee \}.
\]
If $M_\lambda \neq 0$, we say that $\lambda$ is a \emph{weight} of $M$.

A $U_q(\frakg)$-module $M$ is said to be \emph{integrable} if
\begin{enumerate}
    \item $M = \bigoplus_{\lambda \in P} M_\lambda$,
    \item $M_\lambda$ is finite-dimensional for any $\lambda \in P$, and
    \item $e_i$ and $f_i$ act locally nilpotent on $M$ for any $i \in I$.
\end{enumerate}
In this thesis we will only consider integrable modules.

For a dominant weight $\lambda \in P_+$, let $V(\lambda)$ be the \emph{integrable highest weight module} of highest weight $\lambda$; we denote its highest weight vector by $u_\lambda$.
The elements in the orbit $W\lambda$ are called \emph{extremal weights} of $V(\lambda)$.
It is well-known \cite{Lusztig:IntroductionToQuantumGroups} that $V(\lambda)$ is simple and that it is isomorphic to
\[
    \left. U_q(\frakg) \middle/ \middle(\sum\nolimits_{i\in I} U_q(\frakg)e_i + \sum\nolimits_{h \in P^\vee} U_q(\frakg)(q^h - q^{\langle h, \lambda\rangle}) + \sum\nolimits_{i\in I} U_q(\frakg)f_i^{\langle h_i, \lambda\rangle + 1} \right).
\]

There is a unique non-degenerate symmetric bilinear form $(\cdot, \cdot)$ on $V(\lambda)$ with values in $\bbQ$ \cite{Kashiwara:CrystalBases-QuantizedUniversalEnvelopingAlgebras} such that
\begin{enumerate}
    \item $(u_\lambda, u_\lambda) = 1$,
    \item $(q^h u, v) = (u, q^h v)$ for $u, v \in V(\lambda)$ and $h\in P^*$,
    \item $(e_i u, v) = (u, q_i t_i^{-1} f_i v)$ for $u, v \in V(\lambda)$ and $i\in I$,
    \item $(f_i u, v) = (u, q_i^{-1} t_i e_i v)$ for $u, v \in V(\lambda)$ and $i\in I$.
\end{enumerate}

The \emph{dual} $V(\lambda)^\vee$ of $V(\lambda)$ is, as vector space, defined as the graded dual with respect to the weight space decomposition; the $U_q(\frakg)$-module structure is defined via the vector space isomorphism $V(\lambda) \to V(\lambda)^\vee, u \mapsto (u, \cdot)$.
In particular, $V(\lambda)^\vee$ is also an integrable highest weight module of highest weight $\lambda$, and $V(\lambda)^\vee \cong V(\lambda)$.
Note that the structure differs from the usual dual module structure for Hopf algebras, which is defined via the antipode.
Such a module would have highest weight $w_0\lambda$ if $\frakg$ is of finite type, which is not isomorphic to $V(\lambda)$ in general; if $\frakg$ is not of finite type, it would not be a highest weight module.

Let $w \in W$ and fix reduced decomposition $w = s_{i_1}\cdots s_{i_r}$.
We define the vector
\[
    u_{w\lambda} := f_{i_1}^{(a_1)} \cdots f_{i_r}^{(a_r)} u_\lambda,
\]
where $a_k = \langle h_{i_k}, s_{i_{k+1}}\cdots s_{i_r} \lambda\rangle$ for $k \in [1, r]$.
This definition is independent of the choice of the reduced decomposition of $w$ due to the quantum Verma identities \cite[Proposition 39.3.7]{Lusztig:IntroductionToQuantumGroups}.

The \emph{Demazure module} $V_w(\lambda)$ is the $U_q^+(\frakg)$-submodule of $V(\lambda)$ generated by $u_{w\lambda}$:
\[
    V_w(\lambda) = U_q^+(\frakg) u_{w\lambda}.
\]

\subsection{Canonical basis}

We recall the definition of the canonical basis of $U_q^-(\frakg)$ and $V(\lambda)$ following Kashiwara \cite{Kashiwara:CrystalBases-QuantizedUniversalEnvelopingAlgebras}.

\subsubsection*{$\bbZ$-forms}

Let $U_\bbZ(\frakg)$ be the $\bbZ[q, q^{-1}]$-subalgebra of $U_q(\frakg)$ generated by $e_i^{(n)}$, $f_i^{(n)}$, $t_i$, and $\qbinom{t_i; a}{n}$ ($i\in I$, $a\in \bbZ$, $n \ge 0$) .
Let $U_\bbZ^-(\frakg)$ be the $\bbZ[q,q^{-1}]$-subalgebra of $U_q^-(\frakg)$ generated by $f_i^{(n)}$.

For a dominant weight $\lambda$, we define
\begin{equation}
    V_\bbZ(\lambda) := U_\bbZ^-(\frakg) u_\lambda.
\end{equation}

\subsubsection*{$q$-Boson algebra}
\label{section:definition-q-boson-algebra}

For $i\in I$, let $e_i', e_i''\colon U_q^-(\frakg) \to U_q^-(\frakg)$ be the skew-derivations given by
\[
    e_i'f_j = q_i^{-\langle h_i, \alpha_j\rangle} f_j e_i' + \delta_{ij} \text{ and } e_i'' f_j = q_i^{\langle h_i, \alpha_j\rangle} f_j e_i'' + \delta_{ij}.
\]
Here $f_j$ acts on $U_q^-(\frakg)$ by left multiplication.

Kashiwara \cite{Kashiwara:CrystalBases-QuantizedUniversalEnvelopingAlgebras} showed that $e_i'$ and $e_i''$ that are uniquely determined by the equation
\begin{equation}
    \label{equation:e-commutator}
    [e_i, u] = \frac{t_i e_i'' u - t_i^{-1} e_i' u}{q_i - q_i^{-1}} \quad (u \in U_q^-(\frakg)).
\end{equation}

The \emph{$q$-Boson algebra $B_q(\frakg)$} is the $\bbQ(q)$-algebra generated by $e_i'$, $f_i$ ($i\in I$) subject to the relations
\begin{enumerate}
    \item $e_i' f_j = q_i^{-\langle h_i, \alpha_j\rangle} f_j e_i' + \delta_{ij}$,
    \item $\sum_{k=0}^a (-1)^k e_i'^{(k)} e_j' e_i'^{(a-k)} = \sum_{k=0}^a (-1)^k f_i^{(k)} f_j f_i^{(a-k)} = 0$ for $i \neq j$, $a = 1 - \langle h_i, \alpha_j\rangle$.
\end{enumerate}

$U_q^-(\frakg)$ is a left $B_q(\frakg)$-module \cite[Lemma 3.4.2]{Kashiwara:CrystalBases-QuantizedUniversalEnvelopingAlgebras}.

\subsubsection*{Inner product}

By \cite{Kashiwara:CrystalBases-QuantizedUniversalEnvelopingAlgebras} there is a unique non-degenerate symmetric bilinear form $(\cdot, \cdot)$ on $U_q^-(\frakg)$ with values in $\bbQ$ such that
\[
    (1, 1) = 1 \text{ and } (f_iu, v) = (u, e_i'v) \text{ for } u, v\in U_q^-(\frakg).
\]

\subsubsection*{Canonical basis of $U_q^-(\frakg)$}
\label{section:canonical-basis-definition}

Kashiwara \cite{Kashiwara:CrystalBases-QuantizedUniversalEnvelopingAlgebras} showed that for any $i\in I$
\[
    U_q^-(\frakg) = \bigoplus_{n\ge 0} f_i^{(n)} \ker e_i'.
\]
Let $\te_i, \tf_i\colon U_q^-(\frakg) \to U_q^-(\frakg)$ be the operators defined by
\[
    \te_i(f_i^{(n)}u) = f_i^{(n-1)}u
    \text{ and }
    \tf_i(f_i^{(n)}u) = f_i^{(n+1)}u
    \text{ for } u \in \ker e_i'.
\]
They are commonly referred to as \emph{Kashiwara operators}.

Let $L(\infty)$ be the $A$-submodule of $U_q^-(\frakg)$ generated by $\tf_{i_1}\cdots \tf_{i_r} \cdot 1$ for $r \ge 0$ and $i_1, \dots, i_r \in I$.
Denote by $B(\infty) \subseteq L(\infty)/qL(\infty)$ the images of vectors of the form $\tf_{i_1}\cdots \tf_{i_r} \cdot 1$ and write $b_\infty \in B(\infty)$ for the image of $1$ under the projection $L(\infty) \to L(\infty)/qL(\infty)$.

There is an isomorphism \cite[Theorem 6]{Kashiwara:CrystalBases-QuantizedUniversalEnvelopingAlgebras}
\[
    G\colon L(\infty)/qL(\infty) \to U_\bbQ^-(\frakg) \cap L(\infty) \cap \overline{L(\infty)}
\]
of $\bbQ$-vector spaces.

The image of $B(\infty)$ under $G$ is called the \emph{(lower) global crystal basis} or \emph{canonical basis} of $U_q^-(\frakg)$; we denote it by $\calB(\infty)$.
We refer to its dual basis $\calB(\infty)^\vee$ in the graded dual $U_q^-(\frakg)^\vee$ as the \emph{dual canonical basis}.

\begin{theorem}[{\cite[Theorem 2.1.1]{Kashiwara:CrystalBase-Littelmann-DemazureCharacterFormula}}]
    \label{theorem:crystal-base-star-invariance}
    The crystal base $B(\infty)$ is invariant under the $*$-involution, $*B(\infty) = B(\infty)$.
\end{theorem}

\begin{corollary}
    \label{corollary:canonical-basis-star-invariance}
    The canonical basis $\calB(\infty)$ is invariant under the $*$-involution; more precisely, $*G(b) = G(*b)$ for any $b \in B(\infty)$.
\end{corollary}

From Propositions 5.1.2 and 5.1.3 in \cite{Kashiwara:CrystalBases-QuantizedUniversalEnvelopingAlgebras} and the definition of the canonical basis, we obtain the following criterion.
\begin{theorem}
    \label{theorem:canonical-basis-criterion}
    Let $u \in U_q^-(\frakg)$. If $u \in U_\bbZ^-(\frakg)$, $\bar u = u$ and $(u, u) \in 1+qA$, then either $u$ or $-u$ is in $\calB(\infty)$.
\end{theorem}

\begin{theorem}[{\cite{Kashiwara:CrystalBases-QuantizedUniversalEnvelopingAlgebras}}]
    \label{theorem:matrix-coefficients-lower-global-crystal-basis}
    Let $b \in B(\infty)$ and $i \in I$.
    We have
    \begin{equation}
        \label{equation:e-left-multiplication-canonical-basis}
        e_i G(b) = [\varphi_i(b)+1]_i G(\te_i b) + \sum E_{b,b'}^i G(b'),
    \end{equation}
    where $E_{b,b'}^i \in \bbZ[q,q^{-1}]$ is zero unless $\varphi_i(b') > \varphi_i(b)+1$; and
    \begin{equation}
        \label{equation:f-left-multiplication-canonical-basis}
        f_i G(b) = [\epsilon_i(b)+1]_i G(\tf_i b) + \sum F_{b,b'}^i G(b')
    \end{equation}
    where $F_{b,b'}^i \in \bbZ[q,q^{-1}]$ is zero unless $\epsilon_i(b') > \epsilon_i(b)+1$.
\end{theorem}

\subsubsection*{Canonical basis of $V(\lambda)$}
\label{section:canonical-basis-module-definition}

Let $\lambda$ be a dominant weight.
From the theory of $U_q(\fraksl_2)$ it is known that for any $i\in I$ there is a direct sum decomposition
\[
    V(\lambda) = \bigoplus_{n\ge 0} f_i^{(n)} \ker e_i.
\]

Let $\te_i, \tf_i\colon V(\lambda) \to V(\lambda)$ be the operators defined by
\[
    \te_i (f_i^{(n)} u) = f_i^{(n-1)} u \text{ and } \tf_i (f_i^{(n)} u) = f_i^{(n+1)} u \text{ for } u \in \ker e_i.
\]

Let $-\colon V(\lambda) \to V(\lambda)$, be the involution given by $\overline{F u_\lambda} = \bar F u_\lambda$.

Let $L(\lambda)$ be the $A$-submodule of $V(\lambda)$ generated by $\tf_{i_1}\cdots \tf_{i_r} \cdot 1$ for $r \ge 0$ and $i_1, \dots, i_r \in I$.
Denote by $B(\lambda) \subseteq L(\lambda)/qL(\lambda)$ the images of vectors of the form $\tf_{i_1}\cdots \tf_{i_r} \cdot 1$ and write $b_\lambda \in B(\lambda)$ for the image of $u_\lambda$ under the projection $L(\lambda) \to L(\lambda)/qL(\lambda)$.

There is an isomorphism \cite[Theorem 6]{Kashiwara:CrystalBases-QuantizedUniversalEnvelopingAlgebras}
\[
    G_\lambda\colon L(\lambda)/qL(\lambda) \to V_\bbQ(\lambda) \cap L(\lambda) \cap \overline{L(\lambda)}
\]
of $\bbQ$-vector spaces.

The image of $B(\lambda)$ under $G_\lambda$ is called the \emph{global crystal basis} or \emph{canonical basis} of $V(\lambda)$; we denote it by $\calB(\lambda)$.
We refer to its dual basis $\calB(\lambda)^\vee$ in the graded dual $V(\lambda)^\vee$ as the \emph{dual canonical basis} of $V(\lambda)$.

The $U_q^-(\frakg)$-linear homomorphism $U_q^-(\frakg) \to V(\lambda)$ sending $1$ to $u_\lambda$ induces a surjective homomorphism \cite[Theorem 5]{Kashiwara:CrystalBases-QuantizedUniversalEnvelopingAlgebras}
\begin{equation}
    \pr_\lambda\colon L(\infty)/qL(\infty) \to L(\lambda)/qL(\lambda).
\end{equation}

As a consequence of Theorems 5 and 6 of \cite{Kashiwara:CrystalBases-QuantizedUniversalEnvelopingAlgebras} we have the following theorem.
\begin{theorem}
    \label{theorem:canonical-basis-projection}
    Let $\lambda$ be a dominant weight.
    Then, for any $b\in B(\infty)$,
    \[
        G(b)u_\lambda = G_\lambda(\pr_\lambda(b)).
    \]
    In particular,
    \[
        \calB(\lambda) = \{ G(b)u_\lambda \mid b \in B(\infty), \pr_\lambda(b) \neq 0 \}.
    \]
\end{theorem}

From Propositions 5.1.1 and 5.1.3 in \cite{Kashiwara:CrystalBases-QuantizedUniversalEnvelopingAlgebras} and the definition of the canonical basis of $V(\lambda)$, we obtain the following criterion.
\begin{theorem}
    \label{theorem:canonical-basis-criterion-integrable-highest-weight-module}
    Let $u \in V(\lambda)$. If $u \in V_\bbZ(\lambda)$, $\bar u = u$ and $(u, u) \in 1+qA$, then either $u$ or $-u$ is in $\calB(\infty)$.
\end{theorem}

For $w\in W$ we define the set
\begin{equation}
    B_w(\lambda) := \{ \tf_{i_1}^{a_1} \dots \tf_{i_r}^{a_r} b_\lambda \mid s_{i_1}\cdots s_{i_r} \text{ is a reduced decomposition of } w \}.
\end{equation}

\subsection{Crystals}

The definition of a crystal was first given by Kashiwara \cite{Kashiwara:CrystalBase-Littelmann-DemazureCharacterFormula} as abstraction of the crystal bases $B(\lambda)$ and $B(\infty)$.
We recall the definition given in \emph{loc.cit.}.

A \emph{crystal} is a set $B$ together with maps
\[
    \te_i\colon B \to B \sqcup \{ 0 \},\quad \tf_i\colon B \to B \sqcup \{ 0 \}
\]
\[
    \epsilon_i\colon B \to \bbZ \cup \{ -\infty \}, \varphi_i\colon B \to \bbZ \cup \{ -\infty \}, \weight\colon B \to P
\]
such that for any $i\in I$ and $b \in B$:
\begin{enumerate}
    \item $\varphi_i(b) = \epsilon_i(b) + \langle h_i, \weight(b)\rangle$;
    \item if $\te_i b\in B$, then $\tf_i \te_i b = b$, $\weight(\te_i b) = \weight(b) + \alpha_i$, $\epsilon_i(\te_i b) = \epsilon_i(b) - 1$, and $\varphi_i(\te_i b) = \varphi_i(b) + 1$;
    \item if $\tf_i b\in B$, then $\te_i \tf_i b = b$, $\weight(\tf_i b) = \weight(b) - \alpha_i$, $\epsilon_i(\tf_i b) = \epsilon_i(b) + 1$, and $\varphi_i(\tf_i b) = \varphi_i(b) - 1$;
    \item if $\varphi_i(b) = -\infty$, then $\te_i b = \tf_i b = 0$.
\end{enumerate}
Here $0$ is a ghost element, which is not contained in $B$; the addition on $\bbZ \cup \{ \infty \}$ is extended from $\bbZ$ by
\[
    -\infty + x = -\infty \text{ for any } x \in \bbZ \cup \{ -\infty \}.
\]

A crystal $B$ is called \emph{upper normal} if $\epsilon_i(b) = \max\{ n \ge 0 \mid \te_i^n b \neq 0 \}$ for any $b\in B$ and $i\in I$; similarly, it is called \emph{lower normal} if $\varphi_i(b) = \max\{ n \ge 0 \mid \tf_i^n b \neq 0 \}$ for any $b\in B$ and $i\in I$.
It is called \emph{normal} if it is both upper and lower normal.
Note that if $B$ is upper normal (lower normal), then for any $b\in B$ and $i\in I$ there exists $n \in \bbN$ such that $\te_i^n b = 0$ ($\tf_i^n b = 0$).

A \emph{morphism} between two crystals $B_1$ and $B_2$ is a map $\psi\colon B_1 \sqcup \{0\} \to B_2 \sqcup \{0\}$, such that for $b \in B_1$:
\begin{enumerate}
    \item $\psi(0) = 0$;
    \item if $\psi(b)\in B_2$, then $\weight(\psi(b)) = \weight(b)$, $\epsilon_i(\psi(b)) = \epsilon_i(b)$, $\varphi_i(\psi_i(b)) = \varphi_i(b)$;
    \item if $\psi(b), \psi(\te_i b) \in B_2$, then $\te_i \psi(b) = \psi(\te_i b)$;
    \item if $\psi(b), \psi(\tf_i b) \in B_2$, then $\tf_i \psi(b) = \psi(\tf_i b)$.
\end{enumerate}
A morphism is called \emph{strict} if it commutes with all $\tf_i$ and $\te_i$.

For two crystals $B_1$ and $B_2$ their tensor product $B_1 \otimes B_2$ is defined as follows:
\begin{enumerate}
    \item $B_1 \otimes B_2 := B_1 \times B_2$ as a set and elements $(b_1, b_2) \in B_1 \otimes B_2$ are written as $b_1 \otimes b_2$;
    \item $\epsilon_i(b_1 \otimes b_2) = \max(\epsilon_i(b_1), \epsilon_i(b_2) - \langle h_i, \weight(b_1)\rangle$;
    \item $\varphi_i(b_1 \otimes b_2) = \max(\varphi_i(b_2), \varphi_i(b_1) + \langle h_i, \weight(b_2)\rangle)$;
    \item $\weight(b_1 \otimes b_2) = \weight(b_1) + \weight(b_2)$.
\end{enumerate}
The actions of $\te_i$ and $\tf_i$ are given by
\[
	\te_i (b_1\otimes b_2) =
	\begin{cases}
		\te_i b_1 \otimes b_2, & \text{if } \varphi_i(b_1) \ge \epsilon_i(b_2), \\
		b_1\otimes \te_i b_2, & \text{else};
	\end{cases}
\]
and
\[
	\tf_i (b_1\otimes b_2) =
	\begin{cases}
			\tf_i b_1 \otimes b_2, & \text{if } \varphi_i(b_1) > \epsilon_i(b_2), \\
			b_1\otimes \tf_i b_2, & \text{else};
	\end{cases}
\]
with the convention $b_1 \otimes 0 = 0 \otimes b_2 = 0$.

\subsubsection*{Canonical bases of tensor products}
\label{section:tensor-product-canonical-bases}

Let $\lambda$ and $\mu$ be dominant weights.
We want to define the canonical basis of the tensor product $V(\lambda) \otimes V(\mu)$.
For that purpose, we need to introduce a suitable bar involution on $V(\lambda) \otimes V(\mu)$.
Since the action of $U_q(\frakg)$ on $V(\lambda) \otimes V(\mu)$ is given via the coproduct $\Delta$, we first need to extend the bar involution on $U_q(\frakg)$ to $\Delta(U_q(\frakg))$.

There are two natural ways to do this.
Let $u \in U_q(\frakg)$; the first choice is $\overline{\Delta(u)} := \sum \overline{u_{(1)}} \otimes \overline{u_{(2)}}$ and the second choice is $\Delta(\bar u)$.
Both involutions are linked via the \emph{quasi-$\calR$-matrix} $\Theta$, more precisely, $\Delta(\bar u)\Theta = \Theta\overline{\Delta(u)}$.

To define the quasi-$\calR$-matrix we need to consider the completion
\[
    U_q^+(\frakg) \; \widehat{\otimes} \; U_q^-(\frakg) := \bigoplus_{\xi\in Q} \prod_{\xi=\xi'+\xi''} U_q^+(\frakg)_{\xi'} \otimes U_q^-(\frakg)_{\xi''}.
\]

\begin{theorem}[{\cite[Theorem 4.1.2]{Lusztig:IntroductionToQuantumGroups}}]
    \label{theorem:quasi-r-matrix}
    The quasi-$\calR$-matrix $\Theta$ is the unique element in $U_q^+(\frakg) \;\widehat{\otimes}\; U_q^-(\frakg)$ such that $\epsilon(\Theta) = 1$ and $\Delta(\bar u)\Theta = \Theta\overline{\Delta(u)}$.

    Let $B$ be a basis of $U_q^-(\frakg)$ such that $B_\xi := B \cap U_q^-(\frakg)_\xi$ is a basis of $U_q^-(\frakg)_\xi$ for any $\xi \in Q_-$.
    Let $B^\vee= \{ b^\vee \mid b \in B \}$ be the dual basis of $B$ with respect to $(\cdot, \cdot)$.
    Then, $\Theta = \sum_{\xi\in Q_-} \Theta_\xi$, where
    \[
        \Theta_\xi = \prod_{i\in I} (q_i - q_i^{-1})^{\xi_i} \sum_{b\in B_\xi} (b^\vee)^+ \otimes b \in U_q^+(\frakg)_{-\xi} \otimes U_q^-(\frakg)_\xi,
    \]
    for $\xi = -\sum_{i\in I} \xi_i \alpha_i$ and $+\colon U_q^-(\frakg) \to U_q^+(\frakg)$ is $\bbQ(q)$-algebra isomorphism defined by $f_i \mapsto e_i$ for $i\in I$.
\end{theorem}

\begin{remark}
    The element $\Theta_\xi$ looks different from the one in \cite{Lusztig:IntroductionToQuantumGroups} due to the different choice of coproduct and bilinear form.
\end{remark}

\begin{corollary}[{\cite[Corollary 4.1.3]{Lusztig:IntroductionToQuantumGroups}}]
    We have $\Theta \bar\Theta = \bar\Theta \Theta = 1$.
\end{corollary}

\begin{example}
    If $\frakg$ is of type $A_2$, the quasi-$\calR$-matrix is
    \[
        \Theta = \sum_{n\ge 0} q^{n(n-1)/2} \gamma_n e^{(n)} \otimes f^{(n)},
    \]
    where $\gamma_n = \prod_{k=1}^n (q^k-q^{-k})$.
\end{example}

\begin{example}
    \label{example:quasi-r-matrix-finite-type}
    Assume that $\frakg$ is of finite type and fix a reduced decomposition $w_0 = s_{i_1} \cdots s_{i_N}$.
    For $k\in [1, N]$ put
    \[
        \Theta^{[k]} = \sum_{n\ge 0} q^{n(n-1)/2} \gamma_n T_{i_1} \cdots T_{i_{k-1}}(e_{i_k}^{(n)}) \otimes T_{i_1} \cdots T_{i_{k-1}}(f_{i_k}^{(n)}),
    \]
    where $\gamma_n = \prod_{k=1}^n (q^k-q^{-k})$.
    Then, by \cite[Section 8.30]{Jantzen:Lectures-QuantumGroups}, $\Theta = \Theta^{[1]} \cdots \Theta^{[N]}$.
\end{example}

Let $\lambda$ and $\mu$ be dominant weights.
We define the involution $\Psi\colon V(\lambda) \otimes V(\mu) \to V(\lambda) \otimes V(\mu)$ by
\begin{equation}
    \label{equation:involution-tensor-product}
    \Psi(u\otimes v) = \Theta(\bar{u} \otimes \bar{v}).
\end{equation}

Let $\lambda$ and $\mu$ be dominant weights.
We define a preorder $\blacktriangleright$ on $V(\lambda) \otimes V(\mu)$ by $u \otimes v \blacktriangleright u' \otimes v'$ if and only if
\begin{equation}
    \label{equation:definition-preorder-tensor-product}
    \begin{gathered}
        \weight u' \succ \weight u \text{ and } u' \in U_q^+(\frakg) u, \text{ or } \\
        u' \in \bbQ(q) u, \weight v' \succ \weight v \text{ and } v' \in U_q^+(\frakg) v.
    \end{gathered}
\end{equation}

\begin{theorem}[{\cite[Theorem 27.3.2]{Lusztig:IntroductionToQuantumGroups}}]
    \label{theorem:canonical-basis-tensor-product}
    Let $b_1\in B(\lambda)$ and $b_2\in B(\mu)$. We have
    \[
        G(b_1\otimes b_2) = G_\lambda(b_1)\otimes G_\mu(b_2) + \sum_{b_1'\otimes b_2' \in B(\lambda) \otimes B(\mu)} z_{b_1,b_2;b_1',b_2'} G_\lambda(b_1') \otimes G_\mu(b_2'),
    \]
    where $z_{b_1,b_2;b_1',b_2'} \in q\bbZ[q]$ is zero unless $G_\lambda(b_1') \otimes G_\mu(b_2') \blacktriangleleft G_\lambda(b_1)\otimes G_\mu(b_2)$.
\end{theorem}

Let $\phi\colon V(\lambda + \mu) \to V(\lambda) \otimes V(\mu)$ be the unique $U_q(\frakg)$-module homomorphism sending $u_{\lambda+\mu}$ to $u_\lambda \otimes u_\mu$.
The following lemma relates the involution on $V(\mu + \lambda)$ to the involution defined on the image of $\phi$.

\begin{lemma}
    For $u\in V(\lambda) \otimes V(\mu)$ we have $\phi(\bar u) = \Psi(\phi(u))$.
\end{lemma}
\begin{proof}
    Let $F\in U_q^-(\frakg)$, then
    \[
        \bar{F} \cdot (u_\lambda \otimes u_\mu) = \Delta(\bar F)\Theta(u_\lambda \otimes u_\mu) = \Theta \overline{\Delta(F)} (u_\lambda\otimes u_\mu) = \Psi(F \cdot (u_\lambda\otimes u_\mu)).
    \]
\end{proof}

The most important consequence of the above lemma is that the image of the canonical basis of $V(\lambda + \mu)$ under $\phi$ is the canonical basis of the Cartan component of $V(\lambda) \otimes V(\mu)$.
Thus, we can discuss the canonical basis of $V(\lambda + \mu)$ using the explicit description of the canonical basis of $V(\lambda) \otimes V(\mu)$.

\subsection{The path model}
\label{section:path-model}

We recall the definition of the Littelmann path model following \cite{Littelmann:Paths-RootOperators-RepresentationTheory}.

A \emph{rational piecewise linear path} is a piecewise linear, continuous map $\pi\colon [0, 1] \to P_\bbR := P \otimes_\bbZ \bbR$ such that all turning points are rational.
We identify the paths $\pi$ and $\eta$ if there exists a piecewise linear, non-decreasing, surjective map $\phi\colon [0, 1] \to [0, 1]$ such that $\pi = \phi \circ \eta$.
Denote by $\Pi$ the set of all rational piecewise linear paths $\pi$ with $\pi(0) = 0$ and $\pi(1)\in P$.
We say a path of $\Pi$ is \emph{dominant} if its image is contained in the dominant Weyl chamber.
We write $\pi_\lambda$ for the straight-line path $t \mapsto t\lambda$.

For a path $\pi \in \Pi$ and $i\in I$ we write $H_i\colon [0,1] \to \bbR, t \mapsto \langle h_i, \pi(t)\rangle$.
We say that \emph{$H_i(t)$ attains a local minimum at $t = t_0$} if there exists an $\epsilon > 0$ such that $H_i(t) \ge H_i(t_0)$ for $|t - t_0| < \epsilon$ and $H_i(t) > H_i(t_0)$ for $t_0-\epsilon < t < t_0$ or $t_0 < t < t_0+\epsilon$.

For $i\in I$ let $m_i$ be the global minimum attained by $H_i(t)$ and let $l_i, r_i\colon [0,1] \to [0,1]$ be defined as
\begin{align*}
    l_i(t) &:= \min\{ 1, H_i(s) - m_i \mid t \le s \le 1 \}, \text{ and} \\
    r_i(t) &:= \max\{ 0, 1 + m_i - H_i(s) \mid 0 \le s \le t \}.
\end{align*}

The \emph{root operators} $\te_i, \tf_i\colon \Pi \to \Pi \cup \{ 0 \}$ are given by
\begin{gather*}
    (\te_i \pi)(t) = \begin{cases}
        \pi(t) + r_i(t) \alpha_i, & \text{if } r_i(0) = 0, \\
        0, & \text{otherwise};
    \end{cases} \\
    (\tf_i \pi)(t) = \begin{cases}
        \pi(t) - l_i(t) \alpha_i, & \text{if } l_i(1) = 1, \\
        0, & \text{otherwise}.
    \end{cases}
\end{gather*}

We define
\begin{gather}
    \weight \pi := \pi(1), \\
    \epsilon_i(\pi) := - \min \{ H_i(t) \mid t \in [0, 1] \} \cap \bbZ, \\
    \varphi_i(\pi) := \epsilon_i(\pi) + \langle h_i, \weight \pi \rangle.
\end{gather}

For a dominant path $\eta \in \Pi$ ending in $\lambda \in P_+$ the set
\begin{equation}
    B(\eta) := \{ \tf_{i_1} \dots \tf_{i_r} \eta \mid r \ge 0, i_k \in I \text{ for } k \in [1, r] \} \setminus \{ 0 \}
\end{equation}
is a crystal \cite{Littelmann:LittelwoodRichardsonRule-SymmetrizableKacMoodyAlgebras}.
It was shown by Kashiwara \cite{Kashiwara:Similarity-CrystalBases} that $B(\eta)$ is isomorphic to $B(\lambda)$.
For $\pi \in B(\eta)$ we denote the corresponding element in $B(\lambda)$ by $b_\pi$.

The action of the simple reflections $s_i$ on $B(\eta)$ given by
\[
    s_i \pi =
    \begin{cases}
        \tf_i^{\langle h_i, \weight \pi \rangle} \pi, & \text{if } \langle h_i, \weight \pi \rangle \ge 0, \\
        \te_i^{-\langle h_i, \weight \pi \rangle} \pi, & \text{if } \langle h_i, \weight \pi \rangle < 0,
    \end{cases}
\]
extends to an action of the Weyl group $W$ on $B(\eta)$ \cite{Littelmann:Paths-RootOperators-RepresentationTheory}.

\subsubsection*{Lakshmibai-Seshadri paths}
Let $\lambda$ be a dominant weight.
A \emph{Lakshmibai-Seshadri path (LS path) of shape $\lambda$} is a pair $\pi := (\tau_1, \dots \tau_r; t_1, \dots t_r)$ consisting of a sequence of Weyl group elements $\tau_1 > \ldots > \tau_r$ in $W^\lambda$ and sequence of rational numbers $t_0 := 0 < t_1 < \dots < t_r = 1$, fulfilling the following condition:
for each $k \in [1, r-1]$ there exists a maximal chain
\[
    \tau_k = \sigma_1 \xleftarrow{\beta_1} \sigma_2 \xleftarrow{\beta_2} \dots \xleftarrow{\beta_m} \sigma_{m+1} = \tau_{k+1}
\]
such that
\[
    t_k \langle h_{\beta_j}, \sigma_j \lambda \rangle \in \bbZ \text{ for all } j \in [1, m].
\]
For notational convenience, we allow $\tau_k = \tau_{k+1}$ which is understood as
\[
    (\tau_1, \dots, \tau_k, \tau_{k+2}, \dots, \tau_r; t_1, \dots, t_{k-1}, t_{k+1}, \dots, t_r).
\]

The path $\pi$ is identified with the piecewise linear path $\pi(t)$ defined by
\[
    \pi(t) = (t - t_{n-1}) \tau_n\lambda +  \sum_{k=1}^{n-1} (t_k - t_{k-1}) \tau_k\lambda,
\]
for $t \in [t_{n-1}, t_n]$ and $n \in [1, r]$.

The set $\{ \tau_1, \dots, \tau_r \}$ is called the \emph{support} of $\pi$; we denote it by $\supp \pi$.

The initial direction of $\pi$ is defined as $\iota(\pi) := \tau_1$.
Let $\eta$ be a dominant path ending in $\lambda \in P_+$ such that $\eta$ is a concatenation of LS paths.
If $\pi' \in B(\eta)$ and $\pi \cong \pi'$ via the isomorphism $B(\pi_\lambda) \cong B(\eta)$, we put $\iota(\pi') := \iota(\pi)$.

In this thesis we will only consider LS paths and concatenations thereof.

\section{Base change between canonical basis and path basis}
\label{section:base-change}

\subsection{Modified root operators}
\label{section:modified-root-operators}

For a dominant weight $\lambda \in P_+$ we define a modified version of the root operators $\te_i, \tf_i$ which we will denote by $\ce_i, \cf_i$.

Let $\eta$ be a dominant path and $\pi \in B(\eta)$; put $H_i(t) := \langle h_i, \pi(t)\rangle$.
Let $x \in [0, 1]$ be maximal such that $H_i(t)$ is decreasing on $[0, x]$, let $x^-$ be minimal with $H_i(x^-) = H_i(x)$, and let $x^+$ be maximal such that $H_i(t)$ is increasing on $[x, x^+]$.
Let $y^- \in [0, x^-]$ be minimal such that $H_i(y^-) - H_i(x) \le 1$ and $H_i(t)$ is strictly decreasing on $[y^-, x]$.

Since the definition of the modified root operators only depends on local properties of the path $\pi$, we may assume that $\pi$ is an LS path and write
\[
    \pi = (\tau_1, \dots, \tau_r; t_1, \dots, t_r).
\]

Suppose that $x^- = t_k$ and denote by $\pi^-$ the path
\[
    (\tau_1, \dots, \tau_l, s_i\tau_l, \dots, s_i\tau_k, \tau_{k+1}, \dots, \tau_r; t_1, \dots, t_l, y^-, t_{l+1}, \dots, t_r).
\]
Note that $\pi^-$ is not necessarily an LS path.

We define
\[
    \ce_i \pi = \begin{cases}
        \pi^-, & \text{if } H_i(y^-) - H_i(x) = 1; \\
        0, & \text{otherwise}.
    \end{cases}
\]

Let $y^+ \in [x, x^+]$ be maximal such that $H_i(y^+) - H_i(x) \le 1$ and $H_i(t)$ is strictly increasing on $[x, y^+]$.
Assume that $x = t_k$ and denote by $\pi^+$ the path
\[
    (\tau_1, \dots, \tau_k, s_i\tau_{k+1}, \dots, s_i\tau_l, \tau_l, \dots, \tau_r; t_1, \dots, t_{l-1}, y^+, t_l, \dots, t_r).
\]
Note that $\pi^+$ is not necessarily an LS path.

We define
\begin{equation}
    \label{equation:definition-cfi}
    \cf_i \pi = \begin{cases}
        \pi^+, & \text{if } H_i(y^+) - H_i(x) = 1; \\
        0, & \text{otherwise}.
    \end{cases}
\end{equation}

The following figures illustrate the difference between the action of the root operators $\te_i, \tf_i$ and the modified root operators $\ce_i, \cf_i$ on $H_i(t)$.

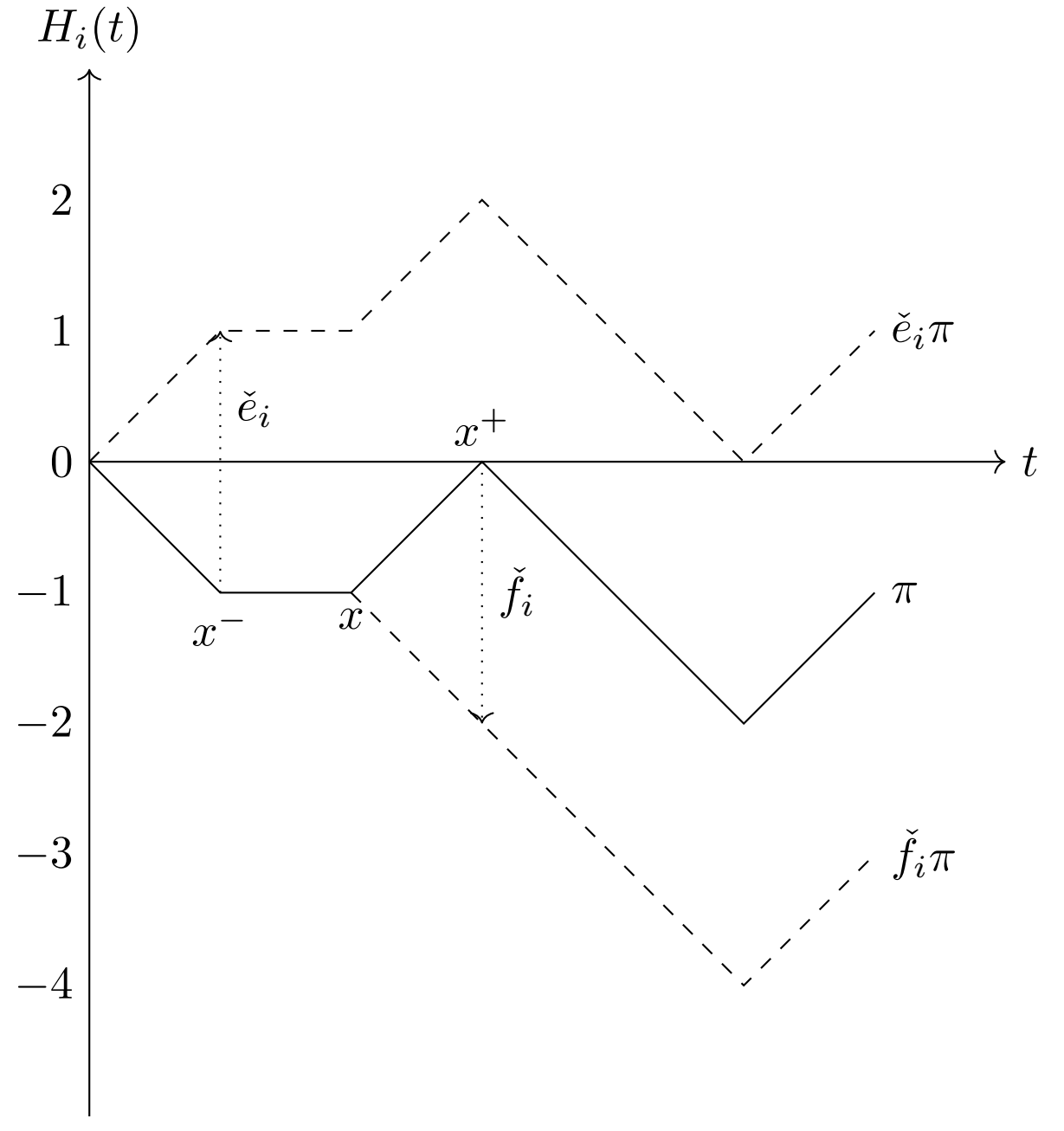
\begin{figure}[H]
    \begin{tikzpicture}
        \draw [->] (0,0) -- (7,0) node[right] {$t$};
        \draw [->] (0,-5) -- (0,3) node[above] {$H_i(t)$};
        \node at (0,2) [left] {$2$};
        \node at (0,1) [left] {$1$};
        \node at (0,0) [left] {$0$};
        \node at (0,-1) [left] {$-1$};
        \node at (0,-2) [left] {$-2$};
        \node at (0,-3) [left] {$-3$};
        \node at (0,-4) [left] {$-4$};

        \node at (1,-1) [below] {$x^-$};
        \node at (2,-1) [below] {$x$};
        \node at (3,0) [above] {$x^+$};
        \draw (0,0) -- (1,-1) -- (2,-1) -- (3,0) -- (5,-2) -- (6,-1) node[right] {$\pi$};
        \draw[dashed] (0,0) -- (1,1) -- (2,1) -- (3,2) -- (5,0) -- (6,1) node[right] {$\ce_i\pi$};
        \draw[->,dotted] (1,-1) -- (1,1) node[pos=0.7, right] {$\ce_i$};
        \draw[dashed] (2,-1) -- (3,-2) -- (5,-4) -- (6,-3) node[right] {$\cf_i\pi$};
        \draw[->,dotted] (3,0) -- (3,-2) node[midway, right] {$\cf_i$};
    \end{tikzpicture}
    \caption{Action of the modified root operators $\ce_i$ and $\cf_i$ on $H_i(t)$.}
\end{figure}

\begin{figure}[H]
    \begin{tikzpicture}
        \draw [->] (0,0) -- (7,0) node[right] {$t$};
        \draw [->] (0,-4) -- (0,2) node[above] {$H_i(t)$};
        \node at (0,1) [left] {$1$};
        \node at (0,0) [left] {$0$};
        \node at (0,-1) [left] {$-1$};
        \node at (0,-2) [left] {$-2$};
        \node at (0,-3) [left] {$-3$};

        \draw (0,0) -- (1,-1) -- (2,-1) -- (3,0) -- (5,-2) -- (6,-1) node[right] {$\pi$};
        \draw[dashed] (4,-1) -- (5,0) -- (6,1) node[right] {$\te_i\pi$};
        \draw[->,dotted] (5,-2) -- (5,0) node[midway, right] {$\te_i$};
        \draw[dashed] (5,-2) -- (6,-3) node[right] {$\tf_i\pi$};
        \draw[->,dotted] (6,-1) -- (6,-3) node[midway, right] {$\tf_{i;1}$};
    \end{tikzpicture}
    \caption{Action of the root operators $\te_i$ and $\tf_i$ on $H_i(t)$.}
\end{figure}
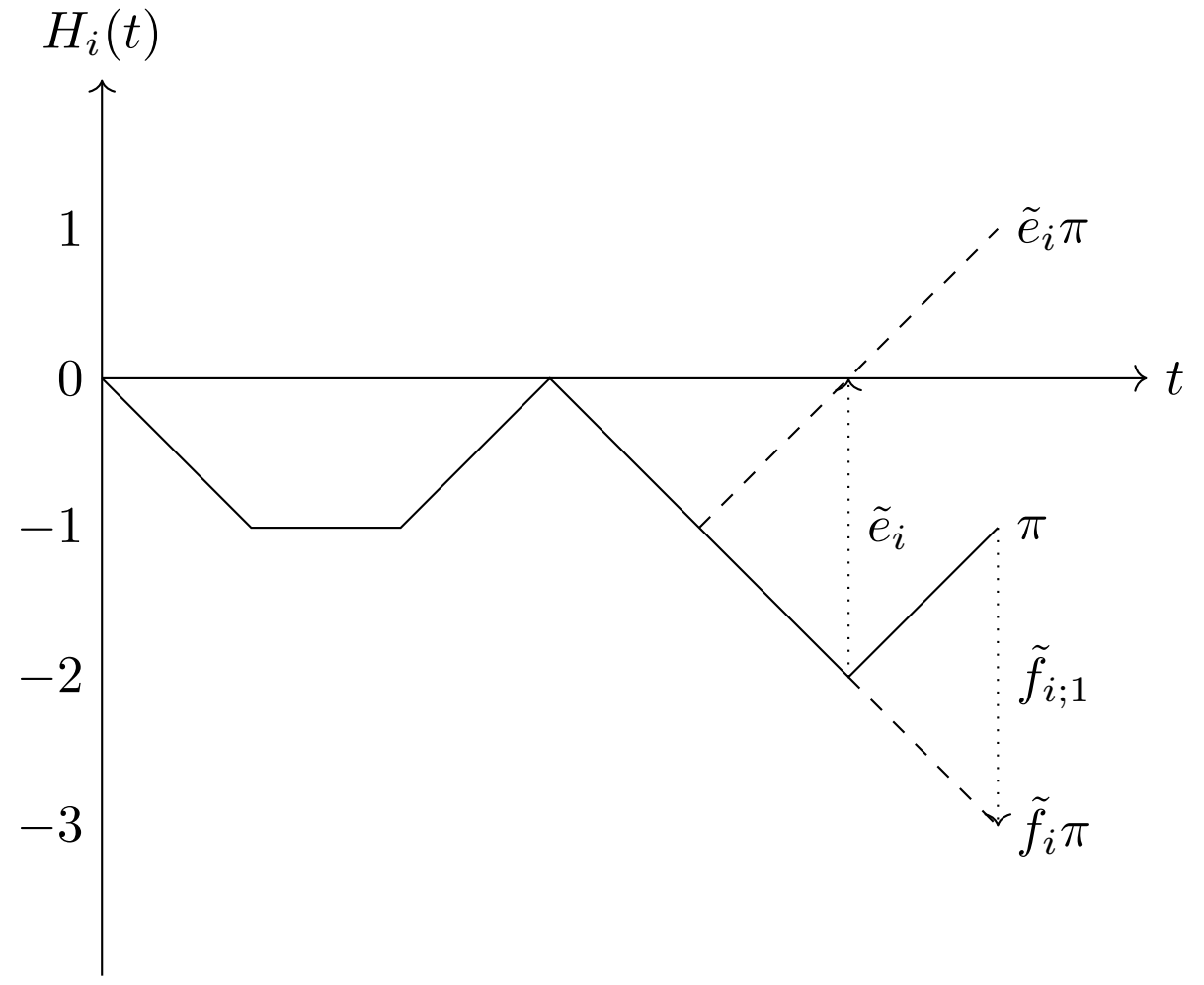

\begin{example}
    Let $\frakg$ be of type $A_2$ and $\lambda = 4\rho$.
    We consider the LS path $\pi = (s_2s_1, s_1, \id; \tfrac38, \tfrac34, 1)$.
    Then,
    \[
        \cf_1\pi = (s_1s_2s_1, s_2s_1, s_1, \id; \tfrac14, \tfrac38, \tfrac34, 1) \text{ and } \tf_1\pi = (s_2s_1, s_1; \tfrac38, 1).
    \]
\end{example}

\begin{lemma}
    \label{lemma:modified-root-operators-well-defined}
    Let $\eta$ be a dominant path.
    The operators $\ce_i, \cf_i$ are well-defined defined maps from $B(\eta)$ to $B(\eta) \sqcup \{ 0 \}$.
\end{lemma}
\begin{proof}
    First, assume that $\eta = \pi_\lambda$ for a dominant weight $\lambda$.
    In this case, the lemma follows directly from the definitions of $\ce_i$ and $\cf_i$ in view of Proposition 4.7 in \cite{Littelmann:Paths-RootOperators-RepresentationTheory}.

    In the general case, $\eta$ is a concatenation of LS paths; here the arguments from the proof of Theorem 7 in \cite{Littelmann:PlacticAlgebra-SemisimpleLieAlgebras} can used to see that the modified root operators are well-defined on $B(\eta)$.
\end{proof}

From the definition of operators $\ce_i, \cf_i$ we get following behaviour for concatenations, that is, tensor products:

Let $\lambda, \mu \in P_+$ be dominant weights, and let $\eta_\lambda$ and $\eta_\mu$ be dominant paths of weight $\lambda$ and $\mu$, respectively.
Consider $\pi \in B(\eta_\lambda)$ and $\pi' \in B(\eta_\mu)$.
If $x$ is the first time that $\langle h_i, (\pi*\pi')(t)\rangle$ attains a local minimum, then
\[
    \ce_i(\pi * \pi') =
    \begin{cases}
        \ce_i\pi * \pi', & \text{if } x \le 1/2; \\
        \pi * \ce_i\pi', & \text{otherwise}.
    \end{cases}
\]
Similarly, we have
\[
    \cf_i(\pi * \pi') =
    \begin{cases}
        \cf_i\pi * \pi', & \text{if } x < 1/2; \\
        \pi * \cf_i\pi', & \text{otherwise}.
    \end{cases}
\]
Here, we use the convention that $0 * \pi' = \pi * 0 = 0$.

From this it is clear that the operators $\ce_i$ and $\cf_i$ are not compatible with the tensor product of crystals.

\begin{lemma}
    We have $\weight(\ce_i\pi) = \weight(\pi) + \alpha_i$ if $\ce_i\pi \neq 0$, and $\weight(\cf_i\pi) = \weight(\pi) - \alpha_i$ if $\cf_i\pi \neq 0$.
\end{lemma}
\begin{proof}
    This follows directly from the definition of $\ce_i$ and $\cf_i$.
\end{proof}

\begin{lemma}
    We have $\cf_i\ce_i\pi = \pi$, if $\ce_i\pi \neq 0$.
\end{lemma}
\begin{proof}
    Let $x^-, x, x^+$, $y^-$ and $y^+$ be as in the definition of $\ce_i$ above; and similarly let $x'^-, x', x'^+$, $y'^-$ and $y'^+$ be the corresponding times for $\pi' := \ce_i\pi$.
    It follows from the definition that $x' = y^-$ and $y'^+ = x^-$.
    Thus, applying $\cf_i$ to $\pi'$ reverts the operation of $\ce_i$ on $\pi$.
\end{proof}

We define
\[
    \cepsilon_i(\pi) := \max\{ k \mid \ce_i^k \pi \ne 0 \}, \quad \cphi_i(\pi) := \max\{ k \mid \cf_i^k \pi \ne 0 \}.
\]

Keep in mind that $B(\eta)$ together with the maps $\ce_i, \cf_i, \cepsilon_i, \cphi_i$ and $\weight$ is not a crystal.

\subsection{Partial orders on the path model}
\label{section:partial-orders}

We recall the partial order $\trianglerighteq$ which was introduced by Littelmann \cite{Littelmann:ContractingModules-StandardMonomialTheory}.
Later, we will see that any linearization of this order works for the base change in Theorem \ref{theorem*:base-change}.
However, for actual proofs there are too many comparable paths.
Therefore, we will introduce a new partial order $\succeq$.

Let $\lambda$ be a dominant weight.
We define all partial orders on $B(\pi_\lambda)$; the definitions generalize naturally to concatentations of LS paths.

The partial order $\triangleright$ is roughly speaking a lexicographic order on LS paths.
For two LS paths $\pi = (\tau_1, \dots, \tau_r; t_1, \dots, t_r)$ and $\pi' = (\sigma_1, \dots, \sigma_s; d_1, \dots, d_s)$ of shape $\lambda$, we write $\pi \triangleright \pi'$, if $\tau_1 > \sigma_1$, or $\tau_1 = \sigma_1$ and $t_1 > d_1$, or $\tau_1 = \sigma_1$, $t_1 = d_1$ and $\tau_2 > \sigma_2$, etc.

\begin{example}
    Assume that $\frakg$ is of type $A_2$ and consider the module $V(2\rho)$.
    We have the paths $\pi_1 = (s_1s_2s_1, \id; \frac12, 1)$, $\pi_2 = (s_2s_1, s_1; \frac12, 1)$, $\pi_3 = (s_2s_1; 1)$, $\pi_4 = (s_1s_2; 1)$. We have $\pi_1 \triangleright \pi_2, \pi_3, \pi_4$; $\pi_3 \triangleright \pi_2$; and the path $\pi_4$ is incomparable to $\pi_2$ and $\pi_3$.
\end{example}

The partial order $\succeq$ on the path model $B(\pi_\lambda)$ is defined as a generalization of the partial order $\succeq$ on the weight lattice.

First, we extend $\succeq$ on the weight lattice $P$ to $P_\bbR = P \otimes_\bbZ \bbR$ in a natural way: For $\lambda, \mu \in P_\bbR$, we say $\lambda \succeq \mu$ if $\lambda - \mu$ is a $\bbR_{\ge0}$-linear combination of simple roots.

Let $\lambda$ be a dominant weight and let $\eta$ be a (arbitrary) dominant path ending in $\lambda$.
For $\pi, \pi' \in B(\eta)$ we write $\pi \succeq \pi'$ if
\begin{equation}
    \label{equation:ls-path-order-succeq}
    \pi(t) - \pi'(t) \succeq 0 \quad \text{for all } t \in [0,1].
\end{equation}

\begin{lemma}
    \label{lemma:partial-order-path-model}
    The relation $\succeq$ is a partial order on $B(\pi_\lambda)$.
\end{lemma}
\begin{proof}
    Reflexivity is clear, for antisymmetry note that $\pi \succeq \pi'$ and $\pi' \succeq \pi$ implies that $\pi(t) - \pi'(t) = 0$ for all $t\in [0, 1]$, and thus $\pi = \pi'$.
    For transitivity, assume that $\pi \succeq \pi'$ and $\pi' \succeq \pi''$.
    Then,
    \[
        \pi(t) - \pi''(t) = (\pi(t) - \pi'(t)) + (\pi'(t) - \pi''(t)) \succeq 0 \quad \text{for all } t \in [0,1].
    \]
\end{proof}

We further define the partial order $\succeq_+$ as
\begin{equation}
    \label{equation:path-order-succeq-plus}
    \pi \succeq_+ \pi' \text{ if and only if } \pi \succeq \pi' \text{ and } \pi \trianglelefteq \pi'.
\end{equation}

\subsection{Bases associated to paths}
\label{section:bases-associated-to-paths}
Littelmann \cite{Littelmann:ContractingModules-StandardMonomialTheory} used the combinatorics of LS paths to construct a monomial basis of $V_\bbZ(\lambda)$.
We recall the construction:

Let $\eta$ be a dominant path ending in $\lambda$ (remember that $\eta$ is assumed to be a concatenation of LS paths).
Let $\pi \in B(\eta)$ and write
\[
    \pi = (\tau_1, \dots; t_1, \dots) * \dots * (\tau_1', \dots; t_1', \dots).
\]
Fix a reduced decomposition $s_{i_1}\dots s_{i_q}$ of $\tau_1$, put $\pi_0 := \pi$ and define
\[
    a_k := \cepsilon_{i_k}(\pi_{k-1}) \text{ and } \pi_k := \ce_{i_k}^{a_k}\pi_{k-1} \quad \text{for } k \in [1, q].
\]
Now, $\pi_q$ has the shape
\[
    (\id; 1) * (\sigma_1, \dots; d_1, \dots) * \dots * (\sigma_1', \dots; d_1', \dots).
\]
Repeating the procedure for the remaining factors, we eventually reach $\pi_r = \eta$.

By \cite{Littelmann:PlacticAlgebra-SemisimpleLieAlgebras}, the sequence $\uli = (i_1, \dots, i_r)$ defines a reduced decomposition of $\iota(\pi)$.
In addition, we have $a_k > 0$ for all $k \in [1, r]$.

We define the \emph{path monomial}
\begin{equation}
    F(\pi, \uli) := f_{i_1}^{(a_1)} \cdots f_{i_r}^{(a_r)}.
\end{equation}

Generally, statements will not depend on specific choice of $\uli$; in such cases we simply write $F(\pi)$.

\begin{example}
    Assume that $\frakg$ is of type $A_2$ and let $\lambda = 4\rho$ and $\eta = \pi_\lambda$.
    For the LS path $\pi := (s_1s_2s_1, s_2s_1, s_1, \id; \tfrac14, \tfrac38, \tfrac34, 1)$ we have $F(\pi, (1,2,1)) = f_1 f_2^{(3)} f_1^{(3)}$ and $F(\pi, (2,1,2)) = f_2^{(2)} f_1^{(4)} f_2$.
\end{example}

For each $\pi \in B(\eta)$ choose a sequence $\uli_\pi$ in $I$ as above.
Then, it was shown by Littelmann \cite{Littelmann:ContractingModules-StandardMonomialTheory} that the set $\{ F(\pi, \uli_\pi) u_\lambda \mid \pi \in B(\eta) \}$ is a basis of $V_\bbZ(\lambda)$.

Since the result is independent of the reduced decomposition chosen for each $\pi \in B(\eta)$, we write
\begin{equation}
    \calP(\eta) := \{ F(\pi) u_\lambda \mid \pi \in B(\eta) \},
\end{equation}
where for each $\pi$ an arbitrary reduced decomposition of $\iota(\pi)$ is fixed.
We call $\calP(\eta)$ the \emph{path monomial basis} associated to $\eta$.

\subsection{Conjectures on the canonical basis}
\label{section:conjectures-on-canonical-basis}

We present the two central conjectures of this paper and show some immediate consequences.

\begin{conjecture}
    \label{conjecture:f-action}
    Let $\lambda$ be a dominant weight and $\eta$ be a dominant path of weight $\lambda$.
    Then, for any path $\pi \in B(\eta)$ and $i \in I$,
    \[
        f_i^{(n)} G_\lambda(b_\pi) = \sum F_{\pi,\pi'}^{i,n} G_\lambda(b_{\pi'}),
    \]
    where $F_{\pi,\pi'}^{i,n} \in \bbZ[q, q^{-1}]$ and the sum ranges over $\pi' \in B(\eta)$ with the following properties:
    \begin{enumerate}
        \item $\pi'(t) \succeq \pi(t) - n\alpha_i$ for all $t \in [0, 1]$.
        \item Let $x \in [0, 1]$ be maximal such that $\pi'(t) - \pi(t)$ is a multiple of $\alpha_i$ for all $t \in [0, x]$.
        Assume that $x < 1$ and let $\epsilon > 0$ be such that on $[x, x+\epsilon]$ the path $\pi$ is supported on $\tau \in W$, and the path $\pi'$ is supported on $\sigma \in W$.
        Then, $\sigma < \max\{\tau, s_i\tau\}$, and if $\pi'(x) - \pi(x) = -n\alpha_i$, then already $\sigma < \tau$.
        \item Let $t_0 \in [0, 1]$ be maximal such that $\langle h_i, \pi(t)\rangle$ is decreasing on $[0, t_0]$, and let $t_1 \in [t_0, 1]$ be maximal such that $\langle h_i, \pi(t)\rangle$ is increasing on $[t_0, t_1]$ and $M := \langle h_i, \pi(t_1)\rangle - \langle h_i, \pi(t_0)\rangle \in \bbZ$.
        If $n \le M$, then $\cf_i^n\pi \neq 0$ and
        \[
            F_{\pi,\cf_i^n\pi}^{i,n} = \qbinom{\cepsilon_i(\pi)+n}{n}_i.
        \]
    \end{enumerate}
\end{conjecture}

\begin{remark}
    Let us make more precise what we mean, when we say that the path $\pi$ is supported on $\tau \in W$ in Conjecture \ref{conjecture:f-action} (ii).
    By assumption $\eta$ is a concatenation of LS paths, that is, $\eta = \pi_{\lambda_1} * \dots * \pi_{\lambda_m}$ for a decomposition of $\lambda$ into dominant weights $\lambda_1$, \dots, $\lambda_m$.
    Thus, $\pi$ can be written as $\pi = \pi_1 * \dots * \pi_m$ with $\pi_k \in B(\pi_{\lambda_k})$ for $k  \in [1, m]$.
    Assume that $\pi$ is parametrized such that $\pi(t) = \pi_1(1) + \dots + \pi_{k-1}(1) + \pi_k(mt - (k-1))$ for $t \in [(k-1)/m, k/m]$.
    Then, if $x < 1$, an $\epsilon > 0$ can be chosen such that $[x, x+\epsilon] \subseteq [(k-1)/m, k/m]$ for some $k \in [1, m]$.
    Assume that $\pi_k = (\tau_1, \dots, \tau_r; t_1 \dots, t_r)$.
    Then, $\epsilon$ can be further refined such that $[mx-(k-1), m(x+\epsilon)-(k-1)] \subseteq [t_{j-1}, t_j]$ for some $j \in [1, r]$.
    We say that on this interval $\pi_k$ is \emph{supported} on $\tau_j \in W^{\lambda_k}$.
    Hence, the support of $\pi$ on $[x, x+\epsilon]$ is $\tau_j$.
\end{remark}

\begin{corollary}
    \label{corollary:f-action-leading-term}
    In Conjecture \ref{conjecture:f-action}, if $\cf_i^n\pi \neq 0$, then $\cf_i^n\pi \trianglerighteq \pi'$ for any $\pi'$ with $F_{\pi,\pi'}^{i,n} \neq 0$.
\end{corollary}
\begin{proof}
    Let $x\in [0, 1]$ be maximal such that $\pi'(t) = \cf_i^n\pi(t)$ for all $t \in [0, x]$ and let $y \in [0, 1]$ be minimal such that $\cf_i^n\pi(y) - \pi(y) = -n\alpha_i$.
    Assume that $\pi' \neq \cf_i^n\pi$.
    Thus, $x < 1$ and there exists $\epsilon > 0$ such that on $[x, x+\epsilon]$ the path $\cf_i^n\pi$ is supported on $\tau \in W$ and the path $\pi'$ is supported on $\sigma \in W$.

    First, consider $x \ge y$.
    Note that $\cf_i^n\pi(t) - \pi(t) = -n\alpha_i$ on $[y, 1]$ and hence $\pi'(x) - \pi(x) = -n\alpha_i$.
    Since $\pi$ is also supported on $\tau$ for $t \in [x, x+\epsilon]$, we have by Conjecture \ref{conjecture:f-action} (ii) that $\sigma < \tau$.
    Hence, $\pi' \triangleleft \cf_i^n\pi$.

    Second, consider $x < y$.
    Note that $\langle h_i, \cf_i^n\pi(t)\rangle$ is decreasing on $[x, x+\epsilon]$, hence $s_i\tau \le \tau$.
    In addition, observe that $\pi$ is supported on $s_i\tau$ or $\tau$.
    Now, Conjecture \ref{conjecture:f-action} (i) implies that we cannot have $\sigma > \tau$.
    Thus, either $\sigma = s_i\tau < \tau$ or it follows from Conjecture \ref{conjecture:f-action} (ii) that $\sigma < \max\{s_i\tau, \tau\} = \tau$.
    Hence, $\pi' \triangleleft \cf_i^n\pi$.
\end{proof}

\begin{example}
    Assume that $\frakg$ is of type $A_2$ and let $\lambda = 4\rho$ and $\eta = \pi_\lambda$.
    \begin{enumerate}
        \item For $\pi = (s_2s_1, s_1, \id; \frac18, \frac14, 1)$ the corresponding element of canonical basis of $V(\lambda)$ is $f_2f_1u_\lambda$.
        We have
        \[
            f_1 (f_2f_1 u_\lambda) = f_1^{(2)}f_2 u_\lambda + f_2 f_1^{(2)} u_\lambda,
        \]
        where the respective paths are
        \[
            \pi_1 := (s_1s_2, \id; \tfrac14, 1) \text{ and } \pi_2 := (s_2s_1, s_1, \id; \tfrac18, \tfrac12, 1).
        \]
        For $\pi_1$ a simple computation gives $\pi_1(\frac18) - \pi(\frac18) = -\frac12\alpha_1 + \frac12\alpha_2$, $\pi_1(\frac14) - \pi(\frac14) = \pi_1(1) - \pi(1) = -\alpha_1$.
        Similarly, we obtain for $\pi_2$ that $\pi_2(\frac18) - \pi(\frac18) = \pi_2(\frac14) - \pi(\frac14) = 0$, and $\pi_2(\frac12) - \pi(\frac12) = \pi_2(1) - \pi(1) = -\alpha_1$.
        Thus, $\pi_1, \pi_2 \succeq \pi - \alpha_1$ and the conjecture holds.
        \item For $\pi = (s_2s_1, s_1, \id; \tfrac38, \tfrac34, 1)$ the corresponding element of canonical basis of $V(\lambda)$ is $f_2^{(3)} f_1^{(3)} u_\lambda$.
        We have
        \[
            f_1 (f_2^{(3)}f_1^{(3)} u_\lambda) = f_2^{(2)} f_1^{(4)} f_2 u_\lambda + f_2^{(3)} f_1^{(4)} u_\lambda.
        \]
        The respective paths are $\pi_1 := (s_1s_2s_1, s_2s_1, s_1, \id; \frac14, \frac38, \frac34, 1)$ and $\pi_2 := (s_1s_2, s_2, \id; \frac12, \frac34, 1)$.
        Here, the conjecture holds, since $\cepsilon_i(\pi) = 0$, $\cf_i\pi = \pi_1$ and $\pi_1 \triangleright \pi_2$.
    \end{enumerate}
\end{example}

The following example shows that, in general, the coefficient of the leading term of Conjecture \ref{conjecture:f-action} is not well-behaved.

\begin{example}
    Assume that $\frakg$ is of type $A_2$ and let $\lambda = 4\rho$ and $\eta = \pi_\lambda$.
    For the path $\pi = (s_1s_2s_1, s_2s_1, s_1, \id; \frac14, \frac12, \tfrac34, 1)$ the corresponding element of the canonical basis is $f_1 f_2^{(4)} f_1^{(3)} u_\lambda$.
    A computation gives
    \[
        f_1^{(2)} (f_1 f_2^{(4)} f_1^{(3)} u_\lambda) = ([4] + [2]) f_2^{(2)} f_1^{(6)} f_2^{(2)} + [3] f_2 f_1^{(6)} f_2^{(3)} u_\lambda + [3] f_2^{(3)} f_1^{(6)} f_2 u_\lambda,
    \]
    where $f_2^{(2)} f_1^{(6)} f_2^{(2)}$ is the element of the canonical basis corresponding to the path $\cf_i^2\pi = (s_1s_2s_1, s_1; \tfrac12, 1)$.
\end{example}

\begin{proposition}
    \label{proposition:base-change-path-monomial-basis-canonical-basis}
    Assume that Conjecture \ref{conjecture:f-action} holds for the dominant path $\eta$ ending in the dominant weight $\lambda$.
    Then, there is a unitriangular base change between the path monomial basis $\calP(\eta)$ and the canonical basis of $V(\lambda)$ for any linearisation of the partial order $\trianglerighteq$ on $B(\eta)$.
\end{proposition}
\begin{proof}
    The proof is by induction on the partial order $\trianglerighteq$.
    The case $\pi = \eta$ is clear, thus we can assume that $\pi \triangleright \eta$ and that the statement holds for all $\pi' \in B(\eta)$ with $\pi' \triangleleft \pi$.
    Fix a sequence $\uli = (i_1, \dots, i_r)$ in $I$ as in the construction of $F(\pi, \uli)$ and put $\pi_k := \cf_{i_k}^{a_k} \cdots \cf_{i_r}^{a_r} \eta$ and $\uli_k = (i_k, \dots, i_r)$ for $k \in [1, r]$.
    We have $\pi = \pi_1 \triangleright \pi_2$, hence using the induction hypothesis we can write
    \[
        F(\pi_2, \uli_2) u_\lambda = G_\lambda(b_{\pi_2}) + \sum_{\pi' \in B(\eta)} z_{\pi_2, \pi'} G_\lambda(b_{\pi'}),
    \]
    where $z_{\pi_2, \pi'} \in \bbZ[q,q^{-1}]$ is zero unless $\pi' \triangleleft \pi_2$.
    We claim that applying $f_{i_1}^{(a_1)}$ to both sides gives
    \[
        F(\pi, \uli) u_\lambda = G_\lambda(b_\pi) + \sum F_{\pi_2,\pi''}^{i_1,a_1} G_\lambda(b_{\pi''}) + \sum z_{\pi_2, \pi'} F_{\pi',\pi'''}^{i_1,a_1} G_\lambda(b_{\pi'''}),
    \]
    where the first sum ranges over $\pi'' \in B(\eta)$ with $\pi'' \triangleleft \pi$ and second sum ranges over $\pi''' \in B(\eta)$ with $\pi''' \triangleleft \pi$.

    Indeed, the coefficient of $G_\lambda(b_\pi)$ is $1$ by Conjecture \ref{conjecture:f-action} (iii), since $a_1 \le \cphi_{i_1}(\pi_2)$ by construction, and $\pi'' \triangleleft \pi$ by Corollary \ref{corollary:f-action-leading-term}.

    To see $\pi''' \triangleleft \pi$, let $x \in [0, 1]$ be maximal with $\pi_2(t) = \pi'(t)$ for all $t \in [0, x]$ and let $y \in [0, 1]$ be maximal such that $\langle h_i, \pi(t)\rangle$ is decreasing on $[0, y]$.
    Note that $\langle h_i, \pi_2(t)\rangle$ is increasing on $[0, y]$.

    If $x \ge y$, then the claim also follows from Corollary \ref{corollary:f-action-leading-term}, since $\pi = \cf_{i_1}^n\pi_2 \triangleright \cf_{i_1}^n\pi'$.
    For $x < y$, let $x' \in [0, 1]$ be maximal such that $\pi'''(t) - \pi'(t)$ is a multiple of $\alpha_{i_1}$ for all $t \in [0, x']$.
    Observe that for $t \in [0, \min\{x, x'\}]$ we have $\pi'''(t) - \cf_i^n\pi(t) \succeq 0$, and if for some $t \in [0, \min\{x, x'\}]$ we have $\pi'''(t) - \cf_{i_1}^n\pi(t) \succ 0$, then $\pi''' \triangleleft \cf_{i_1}^n\pi$.
    Hence, we may assume that $\pi'''(t) = \cf_{i_1}^n\pi(t)$ for all $t \in [0, \min\{x, x'\}]$.

    If $x' < x$, let $\epsilon > 0$ be such that on $[x', x'+\epsilon]$ the path $\pi_2$ is supported on $\tau \in W$, the path $\pi'$ is supported on $\sigma \in W$ and the path $\pi'''$ is supported on $\nu \in W$.
    It follows that $\nu < \max\{\sigma, s_{i_1}\sigma\} = s_{i_1}\tau$, where the inequality is due to Conjecture \ref{conjecture:f-action} (ii).
    Since $\cf_{i_1}^n\pi$ is supported on $s_{i_1}\tau$ for $t\in [x',x'+\epsilon]$, we have $\pi''' \triangleleft \cf_{i_1}^n\pi$.

    If $x' \ge x$, let $\epsilon > 0$ be such that on $[x, x+\epsilon]$ the path $\pi_2$ is supported on $\tau \in W$, the path $\pi'$ is supported on $\sigma \in W$ and the path $\pi'''$ is supported on $\nu \in W$.
    Since $\pi' \triangleleft \pi_2$, we have $\sigma < \tau$ and thus $\nu \le \max\{\sigma, s_{i_1}\sigma\} < s_{i_1}\tau$, where the first inequality is due to Conjecture \ref{conjecture:f-action} (ii).
    Since $\cf_{i_1}^n\pi$ is supported on $s_{i_1}\tau$ for $t\in [x,x+\epsilon]$, we have $\pi''' \triangleleft \cf_{i_1}^n\pi$.
\end{proof}

\begin{conjecture}
    \label{conjecture:e-action}
    Let $\lambda$ be a dominant weight and $\eta$ be a dominant path of weight $\lambda$.
    Then, for any path $\pi \in B(\eta)$ and $i \in I$,
    \[
        e_i G_\lambda(b_\pi) = \sum E_{\pi,\pi'}^i G_\lambda(b_{\pi'})
    \]
    where $E_{\pi,\pi'}^i \in \bbZ[q, q^{-1}]$ and the sum ranges over all $\pi' \in B(\eta)$ with the following properties:
    \begin{enumerate}
        \item $\pi' \succ_+ \pi$.
        \item Let $x \in [0, 1]$ be maximal such that $\pi'(t) - \pi(t)$ is a multiple of $\alpha_j$ for some $j \in I$ and all $t \in [0, x]$.
        Assume that $x < 1$ and let $\epsilon > 0$ be such that on $[x, x+\epsilon]$ the path $\pi$ is supported on $\tau \in W$ and the path $\pi'$ is supported on $\sigma \in W$.
        Then, $\sigma < \max\{\tau, s_j\tau\}$, and if $\pi'(x) - \pi(x) = 0$, then already $\sigma < \tau$.
    \end{enumerate}
\end{conjecture}

\begin{corollary}
    Assume that Conjecture \ref{conjecture:e-action} holds for the dominant path $\eta$ ending in the dominant weight $\lambda$.
    Then, the sets
    \begin{align*}
        V_{\succ_+ \pi}(\lambda) &:= \bigoplus_{\pi'\in B(\eta),  \pi' \succ_+ \pi} \bbZ[q,q^{-1}] G_\lambda(b_{\pi'}) \quad \text{and} \\
        V_{\succeq_+ \pi}(\lambda) &:= \bigoplus_{\pi'\in B(\eta),\pi' \succeq_+ \pi} \bbZ[q,q^{-1}] G_\lambda(b_{\pi'}).
    \end{align*}
    are $U_q^+(\frakg)$-submodules of $V(\lambda)$.
\end{corollary}

Combining Conjectures \ref{conjecture:f-action} and \ref{conjecture:e-action}, we can refine Theorem \ref{theorem:canonical-basis-tensor-product} as follows:
\begin{theorem}
    \label{theorem:canonical-basis-tensor-product-refined}
    Assume that Conjecture \ref{conjecture:e-action} holds for the dominant path $\eta_\lambda$ ending in $\lambda$ and Conjecture \ref{conjecture:f-action} holds for the dominant path $\eta_\mu$ ending in $\mu$.
    Then,
    \[
        G(b_{\pi_1} \otimes b_{\pi_2}) = G_\lambda(b_{\pi_1}) \otimes G_\mu(b_{\pi_2}) + \sum z_{\pi_1, \pi_2;\pi_1', \pi_2'} G_\lambda(b_{\pi_1'}) \otimes G_\mu(b_{\pi_2'}),
    \]
    where $z_{\pi_1, \pi_2;\pi_1', \pi_2'} \in q\bbZ[q]$ is zero unless $\pi_1' \succ_+ \pi_1$ and $\pi_2'(t) \succeq \pi_2(t) - (\weight(\pi_1') - \weight(\pi_1))$ for all $t \in [0, 1]$.
    In particular, $\pi_1'*\pi_2' \succ_+ \pi_1*\pi_2$.
\end{theorem}
\begin{proof}
    Let $\pi_1' \in B(\eta_\lambda)$ and $\pi_2' \in B(\eta_\mu)$ such that $z_{\pi_1, \pi_2;\pi_1', \pi_2'} \neq 0$.
    The fact that $z_{\pi_1, \pi_2;\pi_1', \pi_2'} = 0$, unless $\pi_1' \succ_+ \pi_1$ and $\pi_2'(t) \succeq \pi_2(t) - (\weight(\pi_1') - \weight(\pi_1))$ for all $t \in [0, 1]$, is an immediate consequence of Conjecture \ref{conjecture:e-action} and Conjecture \ref{conjecture:f-action}, respectively.
    Since $\pi_1' \succ_+ \pi_1$, we have for $t \in [0, \frac12]$
    \[
        (\pi_1'*\pi_2')(t) - (\pi_1*\pi_2)(t) = \pi_1'(2t) - \pi_1(2t) \succeq 0;
    \]
    and since $\pi_2' \succ \pi_2 - \xi$, with $\xi := \weight(\pi_1') - \weight(\pi_1)$, we have for $t \in [\frac12, 1]$
    \[
        (\pi_1'*\pi_2')(t) - (\pi_1*\pi_2)(t) = \xi + \pi_2'(2t-1) - \pi_2(2t-1) \succeq \xi + (-\xi) = 0.
    \]
    We conclude that $\pi_1' * \pi_2' \succ_+ \pi_1 * \pi_2$.
\end{proof}

\subsection{Comparsion of bases for fundamental weights}
\label{section:canonical-basis-fundamental-weights}

We investigate the veracity of Conjecture \ref{conjecture:f-action} by comparing the path monomial basis to the canonical basis for fundamental weights.
In this section we assume that $\frakg$ is of finite type.

We first consider fundamental weights of classical type.
A fundamental weight $\varpi_i$ is said to be \emph{of classical type} if $\langle h_\beta, \varpi_i\rangle \le 2$ for all positive roots $\beta$.

\begin{proposition}
    \label{proposition:canonical-basis-classical-fundamental-weight}
    Let $\lambda$ be a fundamental weight of classical type and let $\pi$ be an LS path of shape $\lambda$.
    For any reduced decomposition $\uli = (i_1, \dots, i_r)$ of $\iota(\pi)$ we have
    \[
        G_\lambda(b_\pi) = F(\pi, \uli) u_\lambda.
    \]
\end{proposition}
\begin{proof}
    We prove the statement by induction on the length of $\iota(\pi)$.
    If $\ell(\iota(\pi)) = 0$, then $\pi = \pi_\lambda$ and $G_\lambda(b_\pi) = F(\pi, \uli) u_\lambda = u_\lambda$.
    Now, assume that the length of $\iota(\pi)$ is greater than $0$ and write $F(\pi, \uli) = f_{i_1}^{(a_1)} \cdots f_{i_r}^{(a_r)}$.
    Since $\lambda$ is of classical type, the length of any $\fraksl_2$-chain is at most $2$.
    Therefore, every LS path has the shape $(\tau, \sigma; \frac12, 1)$ with $\tau, \sigma \in W$ (we allow $\tau = \sigma$ and this is understood as the path $(\tau; 1)$).
    From this, one sees that for any $i\in I$ the function $\langle h_i, \pi(t)\rangle$ has a unique minimum if $\epsilon_i(\pi) > 0$; in particular, we have $\epsilon_{i_1}(\pi) = \cepsilon_{i_1}(\pi) = a_1$.
    Using the induction hypothesis on $\pi' := \te_i^{a_1} \pi = \ce_{i_1}^{a_1} \pi$ (note that $\pi' \neq 0$), we obtain
    \[
        G_\lambda(b_{\pi'}) = f_{i_2}^{(a_2)}\cdots f_{i_r}^{(a_r)} u_\lambda.
    \]
    Now, using (\ref{equation:f-left-multiplication-canonical-basis}), we see that
    \[
        f_{i_1}^{(a_1)} \cdots f_{i_r}^{(a_r)} u_\lambda = f_{i_1}^{(a_1)} G_\lambda(b_{\pi'}) = G_\lambda(b_\pi) + \sum F_{\pi,\pi''}^{i_1,a_1} G_\lambda(b_{\pi''}),
    \]
    where $F_{\pi,\pi''}^{i_1,a_1} = 0$ unless $\epsilon_{i_1}(\pi'') > a_1$.
    However, by construction $\weight(\pi') + \alpha_i$ is not a weight of $V(\lambda)$ (the length of any $\fraksl_2$-chain is at most $2$).
    Thus, the only term on the right hand side of the above equation is $G_\lambda(b_\pi)$.
    The proposition follows.
\end{proof}

\begin{proposition}
    \label{proposition:conjecture-classical-fundamental-weight}
    Conjecture \ref{conjecture:f-action} holds for fundamental weights of classical type.
\end{proposition}
\begin{proof}
    Let $\lambda$ be a fundamental weight of classical type and let $\pi = (\tau, \sigma; \frac12, 1)$ be an LS path of shape $\lambda$ (we allow $\tau = \sigma$ for convenience).
    For any $i \in I$ and $n \ge 0$ we have $f_i^{(n)} G_\lambda(b_\pi) = \sum F_{\pi,\pi'}^{i,n} G_\lambda(b_{\pi'}) \in V_w(\lambda)$ by Proposition 3.2.3 of \cite{Kashiwara:CrystalBase-Littelmann-DemazureCharacterFormula}, where $w = \max\{\tau, s_i\tau\}$.
    Thus, if for some $\pi'$ in the above sum $F_{\pi,\pi'}^{i,n} \neq 0$, then $\pi' = (\tau', \sigma'; \frac12, 1)$ for some $\tau' \le w$.
    In particular, $\pi'(\frac12) - \pi(\frac12) = \frac12(\tau'\lambda - \tau\lambda) =: \xi \succeq -\alpha_i$.
    Furthermore, $\pi'(1) - \pi(1) = \xi + \frac12(\sigma'\lambda - \sigma\lambda) = -\alpha_i$.
    This proves (i) and (ii); (iii) follows from Proposition \ref{proposition:canonical-basis-classical-fundamental-weight}, since here $\varphi_i(\pi) = \cphi_i(\pi)$.
    The claim that $F_{\pi,\pi'}^{i,n} \in \bbZ[q, q^{-1}]$ follows directly from Theorem \ref{theorem:matrix-coefficients-lower-global-crystal-basis}.
\end{proof}

\begin{proposition}
    \label{proposition:conjecture-e-action-classical-fundamental-weight}
    Conjecture \ref{conjecture:e-action} holds for fundamental weights of classical type.
\end{proposition}
\begin{proof}
    Let $\lambda$ be a fundamental weight of classical type and let $\pi = (\tau, \sigma; \frac12, 1)$ be an LS path of shape $\lambda$ (we allow $\tau = \sigma$ for convenience).
    For any $i \in I$ we have $e_i G_\lambda(b_\pi) = \sum E_{\pi,\pi'}^i G_\lambda(b_{\pi'}) \in V_\tau(\lambda)$ by Proposition 3.2.3 of \cite{Kashiwara:CrystalBase-Littelmann-DemazureCharacterFormula}.
    Thus, if for some $\pi'$ in the above sum $E_{\pi,\pi'}^i \neq 0$, then $\pi' = (\tau', \sigma'; \frac12, 1)$ for some $\tau' \le \tau$.
    In particular, $\pi'(\frac12) - \pi(\frac12) = \frac12(\tau'\lambda - \tau\lambda) =: \xi \in Q_+$.
    Furthermore, $\pi'(1) - \pi(1) = \xi + \frac12(\sigma'\lambda - \sigma\lambda) = \alpha_i$.
    This proves (i) and (ii).
    The claim that $E_{\pi,\pi'}^i \in \bbZ[q, q^{-1}]$ follows directly from Theorem \ref{theorem:matrix-coefficients-lower-global-crystal-basis}.
\end{proof}

Marsh \cite{Marsh:CanonicalBasis-StandardMonomials} investigated the canonical basis of fundamental weights from the perspective of standard monomials.

\begin{remark}
    Marsh \cite{Marsh:CanonicalBasis-StandardMonomials} used the orientation of the Dynkin diagram of $F_4$ opposite to the one in this paper.
    Also, Marsh's computations were made following Lusztig's conventions \cite{Lusztig:IntroductionToQuantumGroups}, whereas this paper follows Kashiwara's conventions \cite{Kashiwara:CrystalBases-QuantizedUniversalEnvelopingAlgebras}.
    We have adapted Marsh's results accordingly.
\end{remark}

We recall the definition of standard monomials for a fundamental weight of classical type.
Let $\tau, \sigma \in W^\lambda$. The tuple $(\tau, \sigma)$ is called an \emph{admissible pair} if there is a sequence
\[
    \tau = \tau_1 > \dots > \tau_{r+1} = \sigma
\]
such that for each $k \in [1, r]$ there exists $i_k \in I$ with $\tau_k = s_{i_k}\tau_{k+1}$ and $\langle h_{i_k}, \tau_{k+1}\lambda\rangle = 2$.

The \emph{standard monomial} $Q(\tau, \sigma)$ associated to an admissible pair $(\tau, \sigma)$ is defined as follows:
\begin{enumerate}
    \item If $\tau = \sigma$, let $s_{i_1}\cdots s_{i_r}$ be a reduced decomposition of $\tau$ and put $a_k := \langle h_{i_k}, s_{i_{k+1}}\cdots s_{i_r}\lambda\rangle$.
    We define
    \[
        Q(\tau, \tau) := f_{i_1}^{(a_1)} \cdots f_{i_r}^{(a_r)} u_\lambda.
    \]
    The definition is independent of the choice of the reduced decomposition due to the quantum Verma identity.
    \item If $\tau > \sigma$ and $\tau = \tau_1 > \dots > \tau_r = \sigma$ is a sequence as above.
    We define
    \[
        Q(\tau, \sigma) := f_{i_1} \cdots f_{i_r} Q(\sigma, \sigma).
    \]
    The definition is independent of the choice of the sequence.
\end{enumerate}

There is an alternative description of standard monomials \cite[Definition 4.3]{Marsh:CanonicalBasis-StandardMonomials}, which is equivalent to the one given above \cite[Theorem 4.4]{Marsh:CanonicalBasis-StandardMonomials}.
Using this definition it is easy to see that, for the LS path $\pi = (\tau, \sigma; \frac12, 1)$, the standard monomial $Q(\tau, \sigma)$ agrees with the definition of $F(\pi, \uli)$.
In particular, we recover Theorem 5.1 of \cite{Marsh:CanonicalBasis-StandardMonomials}, that is, the standard monomial basis and the canonical basis agree for fundamental weights of classical type.

Let us now consider fundamental weights which are not of classical type.

\begin{example}
    \label{example:g2-fundamental-weights}
    Assume that $\frakg$ is of type $G_2$.
    For both fundament weights the canonical basis agrees with the path basis, as the following explicit description shows.
    The canonical basis of $V({\varpi_1})$ is
    \[
        \{ u_{\varpi_1}, f_1 u_{\varpi_1}, f_2 f_1 u_{\varpi_1}, f_1 f_2 f_1 u_{\varpi_1}, f_1^{(2)} f_2 f_1 u_{\varpi_1}, f_2 f_1^{(2)} f_2 f_1 u_{\varpi_1}, f_1 f_2 f_1^{(2)} f_2 f_1 u_{\varpi_1} \}.
    \]
    The corresponding LS paths of shape $\varpi_1$ are
    \begin{align*}
        \{ (\id; 1), (s_1; 1), (s_2s_1; 1), (s_1s_2s_1, s_2s_1; \tfrac12, 1) \\
        (s_1s_2s_1; 1), (s_2s_1s_2s_1; 1), (s_1s_2s_1s_2s_1; 1) \}.
    \end{align*}
    The canonical basis of $V({\varpi_2})$ is
    \begin{align*}
        \{ & u_{\varpi_2}, f_2 u_{\varpi_2}, f_1 f_2 u_{\varpi_2}, f_1^{(2)} f_2 u_{\varpi_2}, f_1^{(3)} f_2 u_{\varpi_2}, f_2 f_1^{(2)} f_2 u_{\varpi_2}, f_2 f_1^{(3)} f_2 u_{\varpi_2}, \\
        & f_1 f_2 f_1^{(2)} f_2 u_{\varpi_2}, f_2^{(2)} f_1^{(3)} f_2 u_\lambda, f_1 f_2 f_1^{(3)} f_2 u_{\varpi_2}, f_1 f_2^{(2)} f_1^{(3)} f_2 u_{\varpi_2}, \\
        & f_1^{(2)} f_2^{(2)} f_1^{(3)} f_2 u_{\varpi_2}, f_1^{(3)} f_2^{(2)} f_1^{(3)} f_2 u_{\varpi_2}, f_2 f_1^{(3)} f_2^{(2)} f_1^{(3)} f_2 u_{\varpi_2} \}.
    \end{align*}
    The corresponding LS paths of shape $\varpi_2$ are
    \begin{align*}
        \{ (\id; 1), (s_2; 1), (s_1s_2, s_2; \tfrac13, 1), (s_1s_2, s_2; \tfrac23, 1), (s_1s_2; 1), (s_2s_1s_2, s_1s_2, s_2; \tfrac12, \tfrac23, 1), \\
        (s_2s_1s_2, s_1s_2; \tfrac12, 1), (s_1s_2s_1s_2, s_2s_1s_2, s_1s_2, s_2; \tfrac13, \tfrac12, \tfrac23, 1), (s_2s_1s_2; 1), \\
        (s_1s_2s_1s_2, s_2s_1s_2, s_1s_2; \tfrac13, \tfrac12, 1), (s_1s_2s_1s_2, s_2s_1s_2; \tfrac13, 1), \\
        (s_1s_2s_1s_2, s_2s_1s_2; \tfrac23, 1), (s_1s_2s_1s_2; 1), (s_2s_1s_2s_1s_2; 1) \}.
    \end{align*}
\end{example}

\begin{lemma}
    \label{lemma:conjecture-g2-fundamental-weights}
    Conjectures \ref{conjecture:f-action} and \ref{conjecture:e-action} hold for both fundamental weights in type $G_2$.
\end{lemma}
\begin{proof}
    For $\varpi_1$ this follows from Propositions \ref{proposition:conjecture-classical-fundamental-weight} and \ref{proposition:conjecture-e-action-classical-fundamental-weight}, since $\varpi_1$ is of classical type.
    For $\varpi_2$, we present a table for the action of $f_1$ and $f_2$ on the canonical basis of $V(\varpi_2)$.
    The result can then be verified by comparing with Example \ref{example:g2-fundamental-weights}.
    \begin{table}[H]
        \begin{tabular}{l|l|l}
            & $f_1$ & $f_2$ \\
            \hline
            $u_{\varpi_2}$ & 0 & $f_2 u_{\varpi_2}$ \\
            $f_2 u_{\varpi_2}$ & $f_1 f_2 u_{\varpi_2}$ & $0$ \\
            $f_1 f_2 u_{\varpi_2}$ & $[2]_1 f_1^{(2)} f_2 u_{\varpi_2}$ & $0$ \\
            $f_1^{(2)} f_2 u_{\varpi_2}$ & $[3]_1 f_1^{(3)} f_2 u_{\varpi_2}$ & $f_2 f_1^{(2)} f_2 u_{\varpi_2}$ \\
            $f_1^{(3)} f_2 u_{\varpi_2}$ & $0$ & $f_2 f_1^{(3)} f_2 u_{\varpi_2}$ \\
            $f_2 f_1^{(2)} f_2 u_{\varpi_2}$ & $f_1 f_2 f_1^{(2)} f_2 u_{\varpi_2}$ & $0$ \\
            $f_2 f_1^{(3)} f_2 u_{\varpi_2}$ & $f_1 f_2 f_1^{(3)} f_2 u_{\varpi_2}$ & $[2]_2 f_2^{(2)} f_1^{(3)} f_2 u_\lambda$ \\
            $f_1 f_2 f_1^{(2)} f_2 u_{\varpi_2}$ & $[2]_1 f_1 f_2 f_1^{(3)} f_2 u_{\varpi_2}$ & $[3]_1 f_2^{(2)} f_1^{(3)} f_2 u_\lambda$ \\
            $f_2^{(2)} f_1^{(3)} f_2 u_\lambda$ & $f_1 f_2^{(2)} f_1^{(3)} f_2 u_\lambda$ & $0$ \\
            $f_1 f_2 f_1^{(3)} f_2 u_{\varpi_2}$ & $0$ & $f_1 f_2^{(2)} f_1^{(3)} f_2 u_{\varpi_2}$ \\
            $f_1 f_2^{(2)} f_1^{(3)} f_2 u_{\varpi_2}$ & $[2]_1 f_1^{(2)} f_2^{(2)} f_1^{(3)} f_2 u_{\varpi_2}$ & $0$ \\
            $f_1^{(2)} f_2^{(2)} f_1^{(3)} f_2 u_{\varpi_2}$ & $[3]_1 f_1^{(3)} f_2^{(2)} f_1^{(3)} f_2 u_{\varpi_2}$ & $0$ \\
            $f_1^{(3)} f_2^{(2)} f_1^{(3)} f_2 u_{\varpi_2}$ & $0$ & $f_2 f_1^{(3)} f_2^{(2)} f_1^{(3)} f_2 u_{\varpi_2}$ \\
            $f_2 f_1^{(3)} f_2^{(2)} f_1^{(3)} f_2 u_{\varpi_2}$ & $0$ & $0$
        \end{tabular}
    \end{table}
\end{proof}

Marsh \cite{Marsh:CanonicalBasis-StandardMonomials} showed, using a generalized definition for standard monomials, that in $F_4$ for the fundamental weight $\varpi_3$ there exist standard monomials which no longer belong to the canonical basis.
We compare her example to the corresponding element in the path basis.

She considered the standard monomial
\[
    (f_4f_3 - q f_3f_4)f_3f_2f_3f_1f_2f_3 u_\lambda,
\]
which corresponds to the LS path
\[
    \pi = (s_3s_4s_2s_3s_1s_2s_3, s_3s_2s_3s_1s_2s_3; \tfrac12, 1).
\]
Marsh showed that
\[
    ((f_4f_3 - q f_3f_4)f_3f_2f_3f_1f_2f_3 u_\lambda, (f_4f_3 - q f_3f_4)f_3f_2f_3f_1f_2f_3 u_\lambda) = q^2 + q^4.
\]
Hence, by Theorem \ref{theorem:canonical-basis-criterion-integrable-highest-weight-module} is not an element of the canonical basis.

Here, the element of the path monomial basis is $F(\pi, \uli) u_\lambda = f_3^{(2)} f_4 f_2 f_3 f_1 f_2 f_3 u_\lambda$, with $\uli = (3,4,2,3,1,2,3)$.
We compute
\begin{align*}
    & (f_3^{(2)}f_4f_2f_3f_1f_2f_3 u_\lambda, f_3^{(2)}f_4f_2f_3f_1f_2f_3 u_\lambda) \\
    &= (f_4f_2f_3f_1f_2f_3 u_\lambda, f_4f_2f_3f_1f_2f_3 u_\lambda) \\
    &= (1+q^2) (f_2f_3f_1f_2f_3 u_\lambda, f_2f_3f_1f_2f_3 u_\lambda) \\
    &= (1+q^2) (f_3f_1f_2f_3 u_\lambda, f_3f_1f_2f_3 u_\lambda) \\
    &= (1+q^2) (f_1f_2f_3 u_\lambda, f_1f_2f_3 u_\lambda) \\
    &= (1+q^2) (f_2f_3 u_\lambda, f_2f_3 u_\lambda) \\
    &= (1+q^2) (f_3 u_\lambda, f_3 u_\lambda) \\
    &= 1 + q^2.
\end{align*}
Noting that $(2,1,1,1,1,1)$ is also the adapted string of $\pi$ with respect to $\uli$, we conclude that $F(\pi, \uli) = G_{\varpi_3}(b_\pi)$ using Theorem \ref{theorem:canonical-basis-criterion-integrable-highest-weight-module} and Proposition 10.3 of \cite{Littelmann:Cones-Crystals-Patterns}].

In fact, using OSCAR \cite{OSCAR} we were able to verify that in $V(\varpi_3)$, we indeed have $F(\pi, \uli) = G_{\varpi_3}(b_\pi)$ for all $\pi \in B(\pi_{\varpi_3})$ and any choice of $\uli$ such that $s_\uli = \iota(\pi)$.

However, we found an example for $\varpi_2$, where the path basis no longer coincides with the canonical basis.
Consider the LS path
\[
    \pi = (s_3 s_1 s_2 s_3 s_4 s_3 s_1 s_2 s_3 s_1 s_2, s_1 s_2 s_3 s_4 s_3 s_1 s_2 s_3 s_1 s_2, s_2 s_3 s_1 s_2 s_3 s_1 s_2; \tfrac14, \tfrac12, 1)
\]
of weight $-\varpi_1 + 2\varpi_4$ and fix $\uli = (3, 1, 2, 3, 4, 3, 1, 2, 3, 1, 2)$.
Then, associated monomial is
\[
    F(\pi, \uli) = f_3 f_1 f_2^{(2)} f_3^{(2)} f_4^{(2)} f_3 f_1 f_2^{(2)} f_3^{(2)} f_1 f_2.
\]
Before we are able to compute the value of $(F(\pi, \uli) u_\lambda, F(\pi, \uli) u_\lambda)$, we need some intermediate computations.
For this we put
\begin{align*}
    v_1 &= f_1 f_2^{(2)} f_3^{(2)} f_4^{(2)} f_3 f_1 f_2^{(2)} f_3^{(2)} f_1 f_2 u_\lambda, \\
    v_2 &= f_3 f_1 f_2^{(3)} f_3^{(2)} f_4^{(2)} f_1^{(2)} f_2 f_3^{(2)} f_2 u_\lambda.
\end{align*}

First, we compute
\begin{align*}
    & (f_4^{(2)} f_3 f_1 f_2^{(2)} f_3^{(2)} f_1 f_2 u_\lambda, f_4^{(2)} f_3 f_1 f_2^{(2)} f_3^{(2)} f_1 f_2 u_\lambda) \\
    &= (1+q^2+q^4) (f_3 f_1 f_2^{(2)} f_3^{(2)} f_1 f_2, f_3 f_1 f_2^{(2)} f_3^{(2)} f_1 f_2 u_\lambda) \\
    &= (1+2q^2+2q^4+q^6) (f_1 f_2^{(2)} f_3^{(2)} f_1 f_2 u_\lambda, f_4^{(2)} f_1 f_2^{(2)} f_3^{(2)} f_1 f_2 u_\lambda) \\
    &= (1+2q^2+2q^4+q^6) (f_2^{(2)} f_3^{(2)} f_1 f_2 u_\lambda, f_2^{(2)} f_3^{(2)} f_1 f_2 u_\lambda) \\
    &= (1+2q^2+2q^4+q^6) (f_3^{(2)} f_1 f_2 u_\lambda, f_3^{(2)} f_1 f_2 u_\lambda) \\
    &= (1+2q^2+2q^4+q^6) (f_1 f_2 u_\lambda, f_1 f_2 u_\lambda) \\
    &= (1+2q^2+2q^4+q^6) (f_2 u_\lambda, f_2 u_\lambda) \\
    &= 1+2q^2+2q^4+q^6
\end{align*}
and
\begin{align*}
    & (f_4^{(2)} f_3 f_1 f_2^{(2)} f_3^{(2)} f_1 f_2 u_\lambda, f_3 f_4^{(2)} f_1 f_2^{(2)} f_3^{(2)} f_1 f_2 u_\lambda) \\
    &= q^2 (f_3 f_1 f_2^{(2)} f_3^{(2)} f_1 f_2, f_3 f_1 f_2^{(2)} f_3^{(2)} f_1 f_2 u_\lambda) \\
    &= q^2+q^4.
\end{align*}
Next, we compute
\begin{align*}
    & (v_1, v_1) \\
    &= (f_2^{(2)} f_3^{(2)} f_4^{(2)} f_3 f_1 f_2^{(2)} f_3^{(2)} f_1 f_2 u_\lambda, f_2^{(2)} f_3^{(2)} f_4^{(2)} f_3 f_1 f_2^{(2)} f_3^{(2)} f_1 f_2 u_\lambda) \\
    &= (f_3^{(2)} f_4^{(2)} f_3 f_1 f_2^{(2)} f_3^{(2)} f_1 f_2 u_\lambda, f_3^{(2)} f_4^{(2)} f_3 f_1 f_2^{(2)} f_3^{(2)} f_1 f_2 u_\lambda) \\
    &= (f_4^{(2)} f_3 f_1 f_2^{(2)} f_3^{(2)} f_1 f_2 u_\lambda, f_4^{(2)} f_3 f_1 f_2^{(2)} f_3^{(2)} f_1 f_2 u_\lambda) \\
    &+ (q^2 + 2 + q^{-2}) (f_4^{(2)} f_3 f_1 f_2^{(2)} f_3^{(2)} f_1 f_2 u_\lambda, f_3 f_4^{(2)} f_1 f_2^{(2)} f_3^{(2)} f_1 f_2 u_\lambda) \\
    &= (1+2q^2+2q^4+q^6) + (1+3q^2+3q^4+q^6) \\
    &= 2 + 5q^2 + 5q^4 + 2q^6
\end{align*}
and
\begin{align*}
    & (v_1, v_2) \\
    &= (f_2^{(2)} f_3^{(2)} f_4^{(2)} f_3 f_1 f_2^{(2)} f_3^{(2)} f_1 f_2 u_\lambda, f_3 f_2^{(3)} f_3^{(2)} f_4^{(2)} f_1^{(2)} f_2 f_3^{(2)} f_2 u_\lambda) \\
    &= (f_3^{(2)} f_4^{(2)} f_3 f_1 f_2^{(2)} f_3^{(2)} f_1 f_2 u_\lambda, f_3 f_2 f_3^{(2)} f_4^{(2)} f_1^{(2)} f_2 f_3^{(2)} f_2 u_\lambda) \\
    &= (f_4^{(2)} f_3 f_1 f_2^{(2)} f_3^{(2)} f_1 f_2 u_\lambda, f_2 f_3 f_4^{(2)} f_1^{(2)} f_2 f_3^{(2)} f_2 u_\lambda + f_3 f_4^{(2)} f_1 f_2^{(2)} f_3^{(2)} f_1 f_2 u_\lambda) \\
    &= (f_3 f_1 f_2^{(2)} f_3^{(2)} f_1 f_2 u_\lambda, 0 + f_3 f_1 f_2^{(2)} f_3^{(2)} f_1 f_2 u_\lambda) \\
    &= (q^2 + q^4) (f_1 f_2^{(2)} f_3^{(2)} f_1 f_2 u_\lambda, f_1 f_2^{(2)} f_3^{(2)} f_1 f_2 u_\lambda) \\
    &= q^2 + q^4.
\end{align*}
Finally, we compute
\begin{align*}
    & (F(\pi, \uli) u_\lambda, F(\pi, \uli) u_\lambda) \\
    &= (1+q^2) (v_1, v_1 + v_2) \\
    &= (1+q^2) ((2 + 5q^2 + 5q^4 + 2q^6) + (q^2 + q^4)) \\
    &= (1+q^2) (2 + 6q^2 + 6q^4 + 2q^6) \\
    &= 2 + 8q^2 + 12q^4 + 8q^6 + 2q^8.
\end{align*}

Hence, by Theorem \ref{theorem:canonical-basis-criterion-integrable-highest-weight-module}, the vector $F(\pi, \uli) u_\lambda$ is not an element of the canonical basis.

For the path $\pi$ we further verified using OSCAR \cite{OSCAR} that for any choice of $\uli$ such that $s_\uli = \iota(\pi)$, we have $F(\pi, \uli) \ne G_\lambda(b_\pi)$.

\subsection{Proof of the conjecture for standard Young tableaux paths}
\label{section:proof-conjecture-standard-young-tableaux-paths}

In the previous section we verified Conjecture \ref{conjecture:f-action} for fundamental weights of classical type and for both fundamental weights in $G_2$.
We now aim to use tensor products to lift these results to general dominant weights $\lambda = \lambda_1\varpi_1 + \dots + \lambda_N\varpi_N$, where the dominant path is $\eta := \pi_{\lambda_1\varpi_1} * \dots * \pi_{\lambda_N\varpi_N}$.
The elements of $B(\eta)$ are called \emph{standard Young tableaux (paths)} \cite{Littelmann:PlacticAlgebra-SemisimpleLieAlgebras}.

First, let us state some results on the canonical basis of tensor products, which we will need in the proofs of Propositions \ref{proposition:left-action-tensor-product} and \ref{proposition:e-action-tensor-product}.
\begin{lemma}
    \label{lemma:q-power-identity}
    For integers $m, n \ge 0$ we have the following identity:
    \[
        q^{n(m-n)} = \sum_{k=0}^n q^{-k(m+k)} q^{(n-k)(n-k-1)/2} \qbinom{n}{k} \prod_{j=k+1}^n (q^{m+j} - q^{-m-j}).
    \]
\end{lemma}
\begin{proof}
    We prove the lemma by induction on $n$.
    For $n = 0$ the claim is trivial; thus, we may assume that $n > 0$ and that the claim holds for $n-1$.
    Using the $q$-binomial identity
    \[
        \qbinom{n}{k} = q^k \qbinom{n-1}{k} + q^{-n+k} \qbinom{n-1}{k-1},
    \]
    we may split the sum into two parts.
    We will use the notation $\{t\} := q^t - q^{-t}$ for better readability.

    For the first part we compute
    \begin{align*}
        & \sum_{k=0}^{n-1} q^{-k(m+k)} q^{(n-k)(n-k-1)/2} q^k \qbinom{n-1}{k} \prod_{j=k+1}^n \{m+j\} \\
        &= q^{n-1} \{m+n\} \sum_{k=0}^{n-1} q^{-k(m+k)} q^{(n-k-1)(n-k-2)/2} \qbinom{n-1}{k} \prod_{j=k+1}^{n-1} \{m+j\} \\
        &= q^{n-1} \{m+n\} q^{(n-1)(m+n-1)} \\
        &= q^{n(m+n)} - q^{(n-2)(m+n)}.
    \end{align*}

    For the second part we compute
    \begin{align*}
        & \sum_{k=1}^n q^{-k(m+k)} q^{(n-k)(n-k-1)/2} q^{-n+k} \qbinom{n-1}{k-1} \prod_{j=k+1}^n \{m+j\} \\
        &= \sum_{k=0}^{n-1} q^{-(k+1)(m+k+1)} q^{(n-k-1)(n-k-2)/2} q^{-n+k+1} \qbinom{n-1}{k} \prod_{j=k+2}^n \{m+j\} \\
        &= q^{-(m+n)} \sum_{k=0}^{n-1} q^{-k(m+k+1)} q^{(n-k-1)(n-k-2)/2} \qbinom{n-1}{k} \prod_{j=k+1}^{n-1} \{m+j+1\} \\
        &= q^{-(m+n)} q^{(n-1)(m+n)} \\
        &= q^{(n-2)(m+n)}.
    \end{align*}

    Adding both parts together gives the desired identity.
\end{proof}

\begin{proposition}
    \label{proposition:canonical-basis-tensor-product-summands}
    Let $\lambda$ and $\mu$ be dominant weights.
    Let $\psi\colon V(\lambda + \mu) \to V(\lambda) \otimes V(\mu)$ be the embedding given by $u_{\lambda + \mu} \mapsto u_\lambda \otimes u_\mu$, and $\psi_0\colon B(\lambda+\mu) \to B(\lambda) \otimes B(\mu)$ be the embedding given by $b_{\lambda + \mu} \mapsto b_\lambda \otimes b_\mu$.
    Let $b \in B(\lambda+\mu)$.
    Assume that $\psi_0(b) = b_1 \otimes b_2$ and write
    \[
        \psi(G_{\lambda+\mu}(b)) = \sum z_{b_1,b_2;b_1',b_2'} G_\lambda(b_1') \otimes G_\mu(b_2').
    \]
    If $\varphi_i(b_1) > \epsilon_i(b_2)$, then
    \[
        z_{\pi_1,\pi_2;\te_i^n b_1,\tf_i^n b_2} = q_i^{n(\varphi_i(b_1)+n)} \qbinom{\epsilon_i(b_2)+n}{n}_i,
    \]
    and if $\varphi_i(b_1) \le \epsilon_i(b_2)$, then
    \[
        z_{b_1,b_2;\te_i^n b_1,\tf_i^n b_2} = q_i^{n(\epsilon_i(b_2)+n)} \qbinom{\varphi_i(b_1)+n}{n}_i.
    \]
\end{proposition}
\begin{proof}
    Let $\Psi\colon V(\lambda) \otimes V(\mu) \to V(\lambda) \otimes V(\mu)$ be the involution defined in (\ref{equation:involution-tensor-product}).
    The coefficients $z_{b_1,b_2;b_1',b_2'}$ are uniquely determined by the requirements that $\psi(G_{\lambda+\mu}(b))$ is $\Psi$-invariant, $z_{b_1,b_2;b_1,b_2} = 1$, $z_{b_1,b_2;b_1',b_2'} \in q\bbZ[q]$ for $b_1' \otimes b_2' \neq b_1 \otimes b_2$, and $z_{b_1,b_2;b_1',b_2'} = 0$ unless $G_\lambda(b_1') \otimes G_\mu(b_2') \blacktriangleleft G_\lambda(b_1)\otimes G_\mu(b_2)$.
    Therefore, we only need to consider the summands $\Theta_{-n\alpha_i}$ ($n \ge 0$) of the quasi-$\calR$-matrix (cf. Theorem \ref{theorem:quasi-r-matrix}).

    We prove the case $\varphi_i(b_1) > \epsilon_i(b_2)$ by induction on $n$; the case $\varphi_i(b_1) \le \epsilon_i(b_2)$ is similar.
    For $n = 0$ there is nothing to show; thus, let $n > 0$ and assume that the lemma holds for all $k \in [0, n-1]$.
    Note that due to the assumption $\varphi_i(b_1) > \epsilon_i(b_2)$, we have $q_i^{n(\varphi_i(b_1)+n)} \qbinom{\epsilon_i(b_2)+n}{n}_i \in q\bbZ[q]$.

    We claim that $z_{\pi_1,\pi_2;\te_i^n b_1,\tf_i^n b_2} = q_i^{n(\varphi_i(b_1)+n)} \qbinom{\epsilon_i(b_2)+n}{n}_i$.
    Indeed, under this choice the application of $\Psi$ to $\psi(G_{\lambda+\mu}(b))$ gives the summand
    \begin{align*}
        & \sum_{k=0}^n q_i^{-k(\varphi_i(b_1)+k)} q_i^{(n-k)(n-k-1)/2} \qbinom{\epsilon_i(b_2)+n}{n-k}_i \qbinom{\epsilon_i(b_2)+k}{k}_i \\
        & \prod_{j=k+1}^n (q_i^{\varphi_i(b_1)+j} - q_i^{-\varphi_i(b_1)-j}) G_\lambda(b_{\te_i^n b_1}) \otimes G_\mu(b_{\tf_i^n b_2}).
    \end{align*}
    Using the identitiy
    \[
        \qbinom{\epsilon_i(b_2)+n}{n-k}_i \qbinom{\epsilon_i(b_2)+k}{k}_i = \qbinom{\epsilon_i(b_2)+n}{n}_i \qbinom{n}{k}_i.
    \]
    and Lemma \ref{lemma:q-power-identity} the above sum becomes
    \[
        q_i^{n(\varphi_i(b_1)+n)} \qbinom{\epsilon_i(b_2)+n}{n}_i G_\lambda(b_{\te_i^n b_1}) \otimes G_\mu(b_{\tf_i^n b_2}).
    \]
    Now, the uniqueness of $z_{\pi_1,\pi_2;\te_i^n b_1,\tf_i^n b_2}$ implies the claim.
\end{proof}

\begin{lemma}
    \label{lemma:e-action-estimation}
    Let $\lambda$ be a dominant weight and $\eta$ be a dominant path of weight $\lambda$.
    Assume that Conjecture \ref{conjecture:e-action} holds for $V(\lambda)$ and $\eta$.
    Let $\pi \in B(\eta)$ and let $N \subseteq V(\lambda)$ be the minimal $\bbZ[q,q^{-1}]$-submodule such that $U_\bbZ^+(\frakg) G_\lambda(b_\pi) \subseteq N$ and there exists a subset $B' \subseteq B(\eta)$ with $N = \bigoplus_{\pi' \in B'} \bbZ[q,q^{-1}] G_\lambda(b_{\pi'})$.
    Then, any $\pi' \in B'$ has the following property:

    Let $x \in [0, 1]$ be maximal such that $\pi'(t) - \pi(t)$ is a multiple of $\alpha_j$ for some $j \in I$ and all $t \in [0, x]$.
    Assume that $x < 1$ and let $\epsilon > 0$ be such that on $[x, x+\epsilon]$ the path $\pi$ is supported on $\tau \in W$ and the path $\pi'$ is supported on $\sigma \in W$.
    Then, $\sigma < \max\{\tau, s_j\tau\}$, and if $\pi'(x) - \pi(x) = 0$, then already $\sigma < \tau$.
\end{lemma}
\begin{proof}
    Let $\pi'' \in B'$ and put $\xi := \weight(\pi'') - \weight(\pi)$.
    We prove the lemma by induction on $\trace(\xi)$.
    If $\trace(\xi) = 0$, then $\pi'' = \pi$ and $x = 1$, so there is nothing to show.
    Thus, we may assume that $\trace(\xi) > 0$ and that the claim holds for all $\pi' \in B'$ with $\trace(\weight(\pi') - \weight(\pi)) < \trace(\xi)$.
    By the definition of $N$, there must exist $i \in I$ and $\pi' \in B(\eta)_{\weight(\pi)+\xi-\alpha_i} \cap B'$ such that $\xi - \alpha_i \in Q_+$ and $E_{\pi',\pi''}^i \neq 0$.

    Fix $j \in I$.
    Let $x' \in [0, 1]$ be maximal such that $\pi'(t) - \pi(t)$ is a multiple of $\alpha_j$ for all $t \in [0, x']$, and let $x'' \in [0, 1]$ be maximal such that $\pi''(t) - \pi(t)$ is a multiple of $\alpha_j$ for all $t \in [0, x'']$.
    Assume that $x'' < 1$.

    We claim that $x' \ge x''$.
    Indeed, if this was not the case, then we could choose $\epsilon > 0$ such that $x'+\epsilon < x''$ and on $[x', x'+\epsilon]$ the path $\pi$ is supported on $\tau \in W$, the path $\pi'$ is supported on $\sigma \in W$, and the path $\pi''$ is supported on $\nu \in W$.
    Now, Conjecture \ref{conjecture:e-action} (ii) implies that $\nu < \max\{\sigma, s_j\sigma\}$ and the induction hypothesis that $\sigma < \max\{\tau, s_j\tau\}$.
    Hence, we have
    \[
        \max\{\nu, s_j\nu\} < \max\{\sigma, s_j\sigma\} < \max\{\tau, s_j\tau\} = \max\{\nu, s_j\nu\},
    \]
    where the inequalities are strict, since $\sigma \notin \{\nu, s_j\nu\} = \{\tau, s_j\tau\}$; a contradiction.

    Thus, we have that $\pi'(t) - \pi(t)$ is a multiple of $\alpha_j$ for all $t \in [0, x'']$.
    Let $\epsilon > 0$ be such that on $[x'', x''+\epsilon]$ the path $\pi$ is supported on $\tau \in W$ and the path $\pi'$ is supported on $\sigma \in W$, and the path $\pi''$ is supported on $\nu \in W$.
    If $x' > x$, we have that $\sigma \in \{\tau, s_j\tau\}$ and Conjecture \ref{conjecture:e-action} (ii) implies that $\nu < \max\{\sigma, s_j\sigma\} = \max\{\tau, s_j\tau\}$, since $\pi''(t) - \pi'(t)$ is a multiple of $\alpha_j$ for all $t \in [0, x'']$.
    If $x' = x$, then either $\nu \in \{\sigma, s_j\sigma\}$, or we have $\nu < \max\{\sigma, s_j\sigma\}$ by Conjecture \ref{conjecture:e-action} (ii), since $\pi''(t)-\pi'(t)$ is a multiple of $\alpha_j$ for all $t \in [0, x'']$.
    In addition, the induction hypothesis implies that $\sigma < \max\{\tau, s_j\tau\}$.
    Hence, $\nu \le \max\{\sigma, s_j\sigma\} < \max\{\tau, s_j\tau\}$, where the second inequality is strict, since $\sigma \notin \{\tau, s_j\tau\}$.

    If $\pi''(x'') - \pi(x'') = 0$, then we see from
    \[
        0 = \pi''(x'') - \pi(x'') = (\pi''(x'') - \pi'(x'')) + (\pi'(x'') - \pi(x''))
    \]
    that $\pi''(x'') - \pi'(x'') = 0$ and $\pi'(x'') - \pi(x'') = 0$ due to Conjecture \ref{conjecture:e-action} (i).
    If $x' > x''$, we have either $\sigma = \tau$ or $\sigma = s_j\tau < \tau$ by Conjecture \ref{conjecture:e-action} (i) and thus $\nu < \sigma \le \tau$ by Conjecture \ref{conjecture:e-action} (ii).
    If $x' = x''$, we have either $\nu = \sigma$ or $\nu = s_j\sigma < \sigma$, by Conjecture \ref{conjecture:e-action} (i), or $\nu < \sigma$, by Conjecture \ref{conjecture:e-action} (ii).
    Hence, we have $\nu \le \sigma < \tau$ by the induction hypothesis.
\end{proof}

\begin{proposition}
    \label{proposition:left-action-tensor-product}
    Let $\lambda$ and $\mu$ be dominant weights.
    Assume that Conjectures \ref{conjecture:f-action} and \ref{conjecture:e-action} hold for both $V(\lambda)$ and $V(\mu)$ with the respective dominant paths $\eta_\lambda$ and $\eta_\mu$.
    Then, Conjecture \ref{conjecture:f-action} also holds for $V(\lambda + \mu)$ and the path $\eta_\lambda * \eta_\mu$.
\end{proposition}
\begin{proof}
    Let $\psi\colon V(\lambda + \mu) \to V(\lambda) \otimes V(\mu)$ be the embedding given by $u_{\lambda + \mu} \mapsto u_\lambda \otimes u_\mu$ and consider $\pi_1 \in B(\eta_\lambda)$ and $\pi_2 \in B(\eta_\mu)$ such that $\pi_1 * \pi_2 \in B(\eta_\lambda * \eta_\mu)$
    We write
    \[
        \psi(G_{\lambda+\mu}(b_{\pi_1*\pi_2})) = \sum z_{\pi_1,\pi_2;\pi_1',\pi_2'} G_\lambda(b_{\pi_1'}) \otimes G_\mu(b_{\pi_2'}),
    \]
    and
    \[
        \psi(f_i^{(n)} G_{\lambda+\mu}(b_{\pi_1*\pi_2})) = \sum F_{\pi_1',\pi_2';\pi_1'',\pi_2''}^{i,n} z_{\pi_1,\pi_2;\pi_1',\pi_2'} G_\lambda(b_{\pi_1''}) \otimes G_\mu(b_{\pi_2''}).
    \]
    Now, to prove Conjecture \ref{conjecture:f-action} for $V(\lambda + \mu)$ it is sufficient to verify properties (i) - (iii) for $\pi_1''*\pi_2''$ , whenever $z_{\pi_1,\pi_2;\pi_1',\pi_2'} \neq 0$ and $F_{\pi_1,\pi_2;\pi_1'',\pi_2''}^{i,n} \neq 0$.

    \emph{Proof of Conjecture \ref{conjecture:f-action} (i).}
    We need to show that $(\pi_1'' * \pi_2'')(t) \succeq (\pi_1 * \pi_2)(t) - n\alpha_i$ for all $t \in [0, 1]$.
    Consider a fixed $k$ in the sum on the right-hand side of
    \begin{align*}
        f_i^{(n)} (G_\lambda(b_{\pi_1'}) \otimes G_\mu(b_{\pi_2'}))
        &= \sum_{k=0}^n q^{k(n-k)} t_i^k f_i^{(n-k)} G_\lambda(b_{\pi_1'}) \otimes f_i^{(k)} G_\mu(b_{\pi_2'}) \\
        &= \sum_{k=0}^n q^{k(n-k)} F_{\pi_1',\pi_1''}^{i,n-k} F_{\pi_2',\pi_2''}^{i,k} t_i^k G_\lambda(b_{\pi_1''}) \otimes G_\mu(b_{\pi_2''}).
    \end{align*}
    We have by Conjecture \ref{conjecture:f-action} for $V(\lambda)$ and $\eta_\lambda$ that any $\pi_1'' \in B(\eta_\lambda)$ with $F_{\pi_1',\pi_1''}^{i,n-k} \neq 0$ fulfills $\pi_1''(t) - \pi_1'(t) \succeq -(n-k)\alpha_i$ for all $t \in [0, 1]$.
    Similarly, we have by Conjecture \ref{conjecture:f-action} for $V(\mu)$ and $\eta_\mu$ that any $\pi_2'' \in B(\eta_\mu)$ with $F_{\pi_2',\pi_2''}^{i,k} \neq 0$ fulfills $\pi_2''(t) - \pi_2'(t) \succeq -k\alpha_i$ for all $t \in [0, 1]$.
    Thus, we have for $t \in [0, 1/2]$ that
    \begin{align*}
        & (\pi_1''*\pi_2'')(t) - (\pi_1*\pi_2)(t) \\
        &= \pi_1''(2t) - \pi_1(2t) \\
        &= (\pi_1''(2t) - \pi_1'(2t)) + (\pi_1'(2t) - \pi_1(2t)) \\
        &\succeq -(n-k)\alpha_i,
    \end{align*}
    and we have for $t \in [1/2, 1]$ that
    \begin{align*}
        & (\pi_1''*\pi_2'')(t) - (\pi_1*\pi_2)(t) \\
        &= \pi_1''(1) - \pi_1(1) + \pi_2''(2t-1) - \pi_2(2t-1) \\
        &= (\pi_1''(1) - \pi_1'(1)) + (\pi_1'(1) - \pi_1(1)) \\
        &\quad + (\pi_2''(2t-1) - \pi_2'(2t-1)) + (\pi_2'(2t-1) - \pi_2(2t-1)) \\
        &\succeq -(n-k)\alpha_i + (\pi_1'(1) - \pi_1(1)) - k\alpha_i  - (\pi_1'(1) - \pi_1(1)) \\
        &= -n\alpha_i.
    \end{align*}
    Hence, we have $(\pi_1''*\pi_2'')(t) - (\pi_1*\pi_2)(t) \succeq -n\alpha_i$ for all $t \in [0, 1]$.

    \emph{Proof of Conjecture \ref{conjecture:f-action} (ii).}
    We assume $n = 1$; the case $n > 1$ follows by induction.
    Thus, we can write
    \[
        f_i (G_\lambda(b_{\pi_1'}) \otimes G_\mu(b_{\pi_2'})) = \sum F_{\pi_1',\pi_2';\pi_1'',\pi_2''}^{i,1} G_\lambda(b_{\pi_1''}) \otimes G_\mu(b_{\pi_2''}).
    \]
    Fix $\pi_1'' \in B(\eta_\lambda)$ and $\pi_2'' \in B(\eta_\mu)$ with $F_{\pi_1',\pi_2';\pi_1'',\pi_2''}^i \neq 0$.

    Let $x' \in [0, 1]$ be maximal such that $(\pi_1'*\pi_2')(t) - (\pi_1*\pi_2)(t)$ is a multiple of $\alpha_i$ for all $t \in [0, x']$.
    Let $x'' \in [0, 1]$ be maximal such that $(\pi_1''*\pi_2'')(t) - (\pi_1*\pi_2)(t)$ is a multiple of $\alpha_i$ for all $t \in [0, x'']$.
    A similar argument as in the proof of Lemma \ref{lemma:e-action-estimation} shows that $x' \ge x''$.

    Assume that $x'' < 1$ and let $\epsilon > 0$ be such that on $[x'', x''+\epsilon]$ the path $\pi_1*\pi_2$ is supported on $\tau \in W$ and the path $\pi_1'*\pi_2'$ is supported on $\sigma \in W$, and the path $\pi_1''*\pi_2''$ is supported on $\nu \in W$.

    First, consider $x'' < 1/2$.
    Then, our assumption means that on $[2x'', 2x'' + 2\epsilon]$ the path $\pi_1$ is supported on $\tau \in W$ and the path $\pi_1'$ is supported on $\sigma \in W$, and the path $\pi_1''$ is supported on $\nu \in W$.
    If $x' > x''$, we have that $\sigma \in \{\tau, s_i\tau\}$ and Conjecture \ref{conjecture:f-action} (ii) implies that $\nu < \max\{\sigma, s_i\sigma\} = \max\{\tau, s_i\tau\}$, since $\pi_1''(t) - \pi_1'(t)$ is a multiple of $\alpha_i$ for all $t \in [0, x'']$.
    If $x' = x''$, we have that $\sigma < \max\{\tau, s_i\tau\}$ by Lemma \ref{lemma:e-action-estimation}.
    Depending on whether $\pi_1'(t) - \pi_1''(t)$ is still a multiple of $\alpha_i$ or not, we apply Conjecture \ref{conjecture:f-action} (i) or (ii), respectively.
    Hence, $\nu \le \max\{\sigma, s_i\sigma\} < \max\{\tau, s_i\tau\}$, where the second inequality comes from the fact that $\sigma \notin \{\tau, s_i\tau\}$.

    If $(\pi_1''*\pi_2'')(x'') - (\pi_1*\pi_2)(x'') = -\alpha_i$, then we see from
    \begin{align*}
        -\alpha_i
        &= (\pi_1''*\pi_2'')(x'') - (\pi_1*\pi_2)(x'') \\
        &= (\pi_1''(2x'') + \pi_1'(2x'')) + (\pi_1'(2x'') - \pi_1(2x'')).
    \end{align*}
    that $\pi_1''(x'') - \pi_1'(x'') = -\alpha_i$, by Conjecture \ref{conjecture:f-action} (i), and $\pi_1'(x'') - \pi_1(x'') = 0$, by Conjecture \ref{conjecture:e-action} (i).
    Now, similar arguments for $x' > x''$ and $x' = x''$ show $\nu < \tau$.

    Second, consider $x'' \ge 1/2$.
    Then, $\pi_1''(1) - \pi_1(1)$ is a multiple of $\alpha_i$ and hence $\pi_1'(1) - \pi_1(1)$ must be as well.
    In particular, $G_\lambda(b_{\pi_1'}) \otimes G_\mu(b_{\pi_2'})$ is obtained from $G_\lambda(b_{\pi_1}) \otimes G_\mu(b_{\pi_2})$ via the action of
    \[
        \Theta_{-k\alpha_i} = q_i^{k(k-1)/2} \prod_{j=1}^k (q_i^j - q_i^{-j}) e_i^{(k)} \otimes f_i^{(k)}
    \]
    for some $k \ge 0$.
    This follows from the definitions of the quasi-$\calR$-matrix (cf. Theorem \ref{theorem:quasi-r-matrix}), of the involution on $V(\lambda) \otimes V(\mu)$ (cf. (\ref{equation:involution-tensor-product})), and of the canonical basis of $V(\lambda) \otimes V(\mu)$ (cf. Theorem \ref{theorem:canonical-basis-tensor-product}).
    Now, using similar arguments to above we see $\nu < \max\{\tau, s_i\tau\}$.

    If $(\pi_1''*\pi_2'')(x'') - (\pi_1*\pi_2)(x'') = -\alpha_i$, then we see from
    \begin{align*}
        -\alpha_i
        &= (\pi_1''*\pi_2'')(x'') - (\pi_1*\pi_2)(x'') \\
        &= ((\pi_1''*\pi_2'')(x'') - (\pi_1'*\pi_2')(x'')) + ((\pi_1'*\pi_2')(x'') - (\pi_1*\pi_2)(x''))
    \end{align*}
    that $(\pi_1''*\pi_2'')(x'') - (\pi_1'*\pi_2')(x'') = -\alpha_i$, by the proof of (i) above, and $(\pi_1'*\pi_2')(x'') - (\pi_1*\pi_2)(x'') = 0$, by Theorem \ref{theorem:canonical-basis-tensor-product-refined}.
    Hence, by the definition of $\pi_1''*\pi_2''$, either $\pi_1'' = \pi_1'$ and $\pi_2''(2x''-1) - \pi_2'(2x''-1) = -\alpha_i$, or $\pi_2'' = \pi_2'$ and $\pi_1''(1) - \pi_1'(1) = -\alpha_i$.
    Furthermore, $\pi_2'(2x''-1)- \pi_2(2x''-1) = \pi_1(1) - \pi_1'(1) = -k\alpha_i$ (recall that $G_\mu(b_{\pi_2'})$ is obtained from the action of $\Theta_{-k\alpha_i}$).
    In both cases, similar arguments for $x' > x''$ and $x' = x''$ show $\nu < \tau$.

    \emph{Proof of Conjecture \ref{conjecture:f-action} (iii).}
    Let $t_0 \in [0, 1]$ be maximal such that $\langle h_i, (\pi_1*\pi_2)(t)\rangle$ is decreasing on $[0, t_0]$ and let $t_1 \in [t_0, 1]$ be maximal such that $\langle h_i, (\pi_1*\pi_2)(t)\rangle$ is increasing on $[t_0, t_1]$ and $M := \langle h_i, (\pi_1*\pi_2)(t_1)\rangle - \langle h_i, (\pi_1*\pi_2)(t_0)\rangle \in \bbZ$.
    Assume that $n \le M$, then $\tf_{i;1}^n\pi \neq 0$.

    First, consider $t_1 \le 1/2$.
    In this case $\cf_i^n(\pi_1*\pi_2) = \cf_i^n\pi_1 * \pi_2$ and $\cepsilon_i(\pi_1*\pi_2) = \cepsilon_i(\pi_1)$.
    Observe that in $\psi(G_{\lambda+\mu}(b_{\pi_1*\pi_2}))$ the only summand from which the term $G_\lambda(b_{\cf_i^n\pi_1}) \otimes G_\mu(b_{\pi_2})$ can be obtained by the action of $f_i^{(n)}$ is $G_\lambda(b_{\pi_1}) \otimes G_\mu(b_{\pi_2})$.
    We obtain by Conjecture \ref{conjecture:f-action} and Corollary \ref{corollary:f-action-leading-term} that
    \begin{align*}
        & f_i^{(n)} G_\lambda(b_{\pi_1}) \otimes G_\mu(b_{\pi_2}) \\
        &= \qbinom{\cepsilon_i(\pi_1)+n}{n}_i G_\lambda(b_{\cf_i^n\pi_1}) \otimes G_\mu(b_{\pi_2}) + \sum F_{\pi_1;\pi_1'}^{i,n} G_\lambda(b_{\pi_1'}) \otimes G_\mu(b_{\pi_2}),
    \end{align*}
    where the sum ranges over $\pi_1' \in B(\eta_\lambda)$ with $\pi_1' \triangleleft \cf_i^n\pi_1$.

    Second, consider $t_1 > 1/2$.
    In this case $\cf_i^n(\pi_1*\pi_2) = \cf_i^{n_1}\pi_1 * \cf_i^{n_2}\pi_2$ for some $n_1, n_2 \ge 0$ with $n_1+n_2 = n$.
    Observe that $\langle h_i, \pi_1(t)\rangle$ has a unique local minimum; hence, $\ce_i^k\pi_1 = \te_i^k\pi_1$ and $\cf_i^k\pi_1 = \tf_i^k\pi_1$ for all $k \ge 0$.
    Furthermore, note that in $\psi(G_{\lambda+\mu}(b_{\pi_1*\pi_2}))$ the only summands from which the term $G_\lambda(b_{\cf_i^{n_1}\pi_1}) \otimes G_\mu(b_{\cf_i^{n_2}\pi_2})$ can be obtained by the action of $f_i^{(n)}$ are $G_\lambda(b_{\te_i^{n_2-k}\pi_1}) \otimes G_\mu(b_{\cf_i^{n_2-k}\pi_2})$ for $k \in [0, n_2]$.
    This follows from Theorem \ref{theorem:matrix-coefficients-lower-global-crystal-basis} for $G_\lambda(b_{\cf_i^{n_1}\pi_1})$ and from Conjecture \ref{conjecture:f-action} for $G_\mu(b_{\cf_i^{n_2}\pi_2})$.
    For the coefficient of $G_\lambda(b_{\te_i^{n_2-k}\pi_1}) \otimes G_\mu(b_{\cf_i^{n_2-k}\pi_2})$ in $\psi(G_{\lambda+\mu}(b_{\pi_1*\pi_2}))$ we may use an analogue of Proposition \ref{proposition:canonical-basis-tensor-product-summands}, since
    \[
        n_2 \le \langle h_i, \pi_2(2t_1-1)\rangle - \langle h_i, \pi_2(\max\{0, 2t_0-1\})\rangle
    \]
    and thus Conjecture \ref{conjecture:f-action} (iii) is applicable for $f_i^{(n_2-k)} G_\mu(b_{\pi_2})$.

    If $t_0 < 1/2$, then $n_1 = \varphi_i(\pi_1)$, $n_2 = n - n_1$ and $\cepsilon_i(\pi_1*\pi_2) = \cepsilon_i(\pi_1) = \epsilon_i(\pi_1)$.
    Based on our previous considerations, we only need to consider the sum
    \[
        \sum_{k=0}^n q_i^{(n_2-k)(\varphi_i(\pi_1)+n_2-k)} q_i^{k(n-k)} t_i^k f_i^{(n-k)} G_\lambda(b_{\te_i^{n_2-k}\pi_1}) \otimes f_i^{(k)} G_\mu(b_{\cf_i^{n_2-k}\pi_2})
    \]
    in $f_i^{(n)} \psi(G_{\lambda+\mu}(b_{\pi_1*\pi_2}))$.
    Hence, the coefficient of $G_\lambda(b_{\cf_i^{n_1}\pi_1}) \otimes G_\mu(b_{\cf_i^{n_2}\pi_2})$ is
    \begin{align*}
        & \sum_{k=0}^n q_i^{(n_2-k)(n-k)} q_i^{k(n-k)} q_i^{-k(\cepsilon_i(\pi_1)+n_1)} \qbinom{\cepsilon_i(\pi_1)+n_1}{n-k}_i \qbinom{n_2}{k}_i \\
        &= \sum_{k=0}^n q_i^{n_2(n-k) - k(\cepsilon_i(\pi_1)+n_1)} \qbinom{\cepsilon_i(\pi_1)+n_1}{n-k}_i \qbinom{n_2}{k}_i \\
        &= \qbinom{\cepsilon_i(\pi_1)+n}{n}_i.
    \end{align*}
    The last equality follows binomial identities in \cite{Lusztig:IntroductionToQuantumGroups}.

    If $t_0 \ge 1/2$, then $n_1 = 0$, $n_2 = n$ and $\cepsilon_i(\pi_1*\pi_2) = \cepsilon_i(\pi_1) + \cepsilon_i(\pi_2)$.
    Based on our previous considerations, we only need to consider the sum
    \[
        \sum_{k=0}^n q_i^{(n-k)(\cepsilon_i(\pi_2)+n-k)} q_i^{k(n-k)} t_i^k f_i^{(n-k)} G_\lambda(b_{\te_i^{n-k}\pi_1}) \otimes f_i^{(k)} G_\mu(b_{\cf_i^{n-k}\pi_2})
    \]
    in $f_i^{(n)} \psi(G_{\lambda+\mu}(b_{\pi_1*\pi_2}))$.
    Hence, the coefficient of $G_\lambda(b_{\pi_1}) \otimes G_\mu(b_{\cf_i^n\pi_2})$ is
    \begin{align*}
        & \sum_{k=0}^n q_i^{(n-k)(\cepsilon_i(\pi_2)+n-k)} q_i^{k(n-k)} q_i^{-k\epsilon_i(\pi_1)} \qbinom{\cepsilon_i(\pi_1)}{n-k}_i \qbinom{\cepsilon_i(\pi_2)+n}{k}_i \\
        &= \sum_{k=0}^n q_i^{(n-k)(\cepsilon_i(\pi_2)+n) - k\epsilon_i(\pi_1)} \qbinom{\cepsilon_i(\pi_1)}{n-k}_i \qbinom{\cepsilon_i(\pi_2)+n}{k}_i \\
        &= \qbinom{\cepsilon_i(\pi_1)+\cepsilon_i(\pi_2)+n}{n}_i.
    \end{align*}
    The last equality follows binomial identities in \cite{Lusztig:IntroductionToQuantumGroups}.
\end{proof}

\begin{proposition}
    \label{proposition:e-action-tensor-product}
    Let $\lambda$ and $\mu$ be dominant weights.
    Assume that Conjectures \ref{conjecture:f-action} and \ref{conjecture:e-action} hold for both $V(\lambda)$ and $V(\mu)$ with the respective dominant paths $\eta_\lambda$ and $\eta_\mu$.
    Then, Conjecture \ref{conjecture:e-action} also holds for $V(\lambda + \mu)$ and the path $\eta_\lambda * \eta_\mu$.
\end{proposition}
\begin{proof}
    Let $\pi_1 \in B(\eta_\lambda)$ and $\pi_2 \in B(\eta_\mu)$ such that $\pi_1 * \pi_2 \in B(\eta_\lambda * \eta_\mu)$.
    Let $\psi\colon V(\lambda + \mu) \to V(\lambda) \otimes V(\mu)$ be the embedding given by $u_{\lambda + \mu} \mapsto u_\lambda \otimes u_\mu$.
    We write
    \[
        \psi(G_{\lambda+\mu}(b_{\pi_1*\pi_2})) = \sum z_{\pi_1,\pi_2;\pi_1',\pi_2'} G_\lambda(b_{\pi_1'}) \otimes G_\mu(b_{\pi_2'}).
    \]
    Fix $\pi_1' \in B(\eta_\lambda)$ and $\pi_2' \in B(\eta_\mu)$ with $z_{\pi_1,\pi_2;\pi_1',\pi_2'} \neq 0$.
    Since Conjectures \ref{conjecture:f-action} and \ref{conjecture:e-action} hold for both $V(\lambda)$ and $V(\mu)$ with the respective dominant paths $\eta_\lambda$ and $\eta_\mu$, we have by Theorem \ref{theorem:canonical-basis-tensor-product-refined} that $\pi_1'*\pi_2' \succeq_+ \pi_1*\pi_2$.

    \emph{Proof of (i)}.
    We have by Conjecture \ref{conjecture:e-action} (i) that
    \begin{align*}
        e_i (G_\lambda(b_{\pi_1'}) \otimes G_\mu(b_{\pi_2'}))
        &= e_i G_\lambda(b_{\pi_1'}) \otimes t_i^{-1} G_\mu(b_{\pi_2'}) + G_\lambda(b_{\pi_1'}) \otimes e_i G_\mu(b_{\pi_1'}) \\
        &= \sum E_{\pi_1',\pi_2';\pi_1'',\pi_2''}^i G_\lambda(b_{\pi_1''}) \otimes G_\mu(b_{\pi_2''}),
    \end{align*}
    where $E_{\pi_1,\pi_2;\pi_1',\pi_2'}^i \in \bbZ[q, q^{-1}]$ and the sum ranges over all $\pi_1'' \in B(\eta_\lambda)$ and $\pi_2'' \in B(\eta_\mu)$ with $\pi_1''*\pi_2'' \succ_+ \pi_1'*\pi_2'$.
    Hence, $\pi_1''*\pi_2'' \succ_+ \pi_1*\pi_2$ for such $\pi_1''$ and $\pi_2''$.

    \emph{Proof of (ii)}.
    Fix $j \in I$, and $\pi_1'' \in B(\eta_\lambda)$ and $\pi_2'' \in B(\eta_\mu)$ with $E_{\pi_1',\pi_2';\pi_1'',\pi_2''}^i \neq 0$.

    Let $x' \in [0, 1]$ be maximal such that $(\pi_1'*\pi_2')(t) - (\pi_1*\pi_2)(t)$ is a multiple of $\alpha_j$ for all $t \in [0, x']$.
    Let $x'' \in [0, 1]$ be maximal such that $(\pi_1''*\pi_2'')(t) - (\pi_1*\pi_2)(t)$ is a multiple of $\alpha_j$ for all $t \in [0, x'']$.
    A similar argument as in the proof of Lemma \ref{lemma:e-action-estimation} shows that $x' \ge x''$.

    Assume that $x'' < 1$ and let $\epsilon > 0$ be such that on $[x'', x''+\epsilon]$ the path $\pi_1*\pi_2$ is supported on $\tau \in W$ and the path $\pi_1'*\pi_2'$ is supported on $\sigma \in W$, and the path $\pi_1''*\pi_2''$ is supported on $\nu \in W$.

    First, consider $x'' < 1/2$.
    Here $\nu < \max\{\tau, s_j\tau\}$ is an immediate consequence of Lemma \ref{lemma:e-action-estimation} for $V(\lambda)$ and $\eta_\lambda$.
    If $(\pi_1''*\pi_2'')(x'') - (\pi_1*\pi_2)(x'') = 0$, then
    \[
        \pi_1''(2x'') - \pi_1(2x'') = (\pi_1''*\pi_2'')(x'') - (\pi_1*\pi_2)(x'') = 0.
    \]
    Hence $\nu < \tau$ by Lemma \ref{lemma:e-action-estimation} for $V(\lambda)$ and $\eta_\lambda$.

    Second, consider $x'' \ge 1/2$.
    Then, $\pi_1''(1) - \pi_1(1)$ is a multiple of $\alpha_j$ and hence $\pi_1'(1) - \pi_1(1)$ must be as well.
    In particular, $G_\lambda(b_{\pi_1'}) \otimes G_\mu(b_{\pi_2'})$ is obtained from $G_\lambda(b_{\pi_1}) \otimes G_\mu(b_{\pi_2})$ via the action of $\Theta_{-k\alpha_j}$ for some $k \ge 0$.
    If $x' > x''$, we have that $\sigma \in \{\tau, s_j\tau\}$ and Conjecture \ref{conjecture:e-action} (ii) for $V(\mu)$ and $\eta_\mu$ implies that $\nu < \max\{\sigma, s_j\sigma\} = \max\{\tau, s_j\tau\}$, since $\pi_2''(t) - \pi_2'(t)$ is a multiple of $\alpha_j$ for all $t \in [0, 2x''-1]$.
    If $x' = x''$, we have that $\sigma < \max\{\tau, s_j\tau\}$ by Conjecture \ref{conjecture:f-action} (ii) for $V(\mu)$ and $\eta_\mu$ and Conjecture \ref{conjecture:e-action} (ii) for $V(\mu)$ and $\eta_\mu$ implies that $\nu \le \max\{\sigma, s_j\sigma\} < \max\{\tau, s_j\tau\}$, where the second inequality comes from the fact that $\sigma \notin \{\tau, s_j\tau\}$.

    If $(\pi_1''*\pi_2'')(x'') - (\pi_1*\pi_2)(x'') = 0$, then we see from
    \begin{align*}
        0
        &= (\pi_1''*\pi_2'')(x'') - (\pi_1*\pi_2)(x'') \\
        &= ((\pi_1''*\pi_2'')(x'') - (\pi_1'*\pi_2')(x'')) + ((\pi_1'*\pi_2')(x'') - (\pi_1*\pi_2)(x''))
    \end{align*}
    that $(\pi_1''*\pi_2'')(x'') - (\pi_1'*\pi_2')(x'') = 0$, by the proof of (i) above, and $(\pi_1'*\pi_2')(x'') - (\pi_1*\pi_2)(x'') = 0$, by Theorem \ref{theorem:canonical-basis-tensor-product-refined}.
    By the definition of $\pi_1''*\pi_2''$, either $\pi_1'' = \pi_1'$ or $\pi_2'' = \pi_2'$.
    In the latter case, we would have $\pi_1''(1) - \pi_1'(1) = \alpha_j$, which is not possible under our current assumptions.
    Hence, $\pi_1'' = \pi_1'$ and $\pi_2''(2x''-1) - \pi_2'(2x''-1) = 0$.
    Furthermore, $\pi_2'(2x''-1)- \pi_2(2x''-1) = \pi_1(1) - \pi_1'(1) = -k\alpha_j$ (recall that $G_\mu(b_{\pi_2'})$ is obtained from the action of $\Theta_{-k\alpha_j}$).
    Now, similar arguments for $x' > x''$ and $x' = x''$ show $\nu < \tau$.
\end{proof}

A special case of Propositions \ref{proposition:left-action-tensor-product} and \ref{proposition:e-action-tensor-product} is for the tensor product $V(m\lambda) \otimes V(n\lambda)$.
In this case, it is possible to give an explicit description of the embedding $V((m+n)\lambda) \hookrightarrow V(m\lambda) \otimes V(n\lambda)$ in terms of LS paths.

\begin{lemma}
    \label{lemma:crystal-cartan-embedding}
    Let $\lambda$ be a dominant weight and $m, n \ge 1$.
    The (unique) embedding of crystals
    \[
        B((m+n)\lambda) \hookrightarrow B(m\lambda) \otimes B(n\lambda)
    \]
    sending $b_{(m+n)\lambda}$ to $b_{m\lambda} \otimes b_{n\lambda}$ is in the language of LS paths given by
    \[
        (\tau_1, \dots, \tau_r; t_1, \dots, t_r) \mapsto (\tau_1, \dots, \tau_r; t_1', \dots, t_r') \otimes (\tau_1, \dots, \tau_r; t_1'', \dots, t_r'')
    \]
    where $t_k' = \min\{(m+n)t_k/m, 1\}$ and $t_k'' = ((m+n)t_k - m t_k')/n$.
\end{lemma}
\begin{proof}
    We verify that $\pi' := (\tau_1, \dots, \tau_r; t_1', \dots, t_r')$ and $\pi'' := (\tau_1, \dots, \tau_r; t_1'', \dots, t_r'')$ are LS paths of shape $m\lambda$ and $n\lambda$, respectively.

    Let $\tau_1 = w_1 > \dots > w_{p+1} = \id$ be a maximal chain on which $\pi$ is supported, and let $\beta_1, \dots, \beta_p$ be positive roots such that $w_{k+1} = \beta_k w_k$ for $k \in [1, p]$.
    Assume $t_1 < m/(m+n)$, then we have for $\pi'$:
    \[
        t_1'\langle h_{\beta_1}, \tau_1 m\lambda\rangle  = \frac{(m+n)t_1}{m}\langle h_{\beta_1}, \tau_1 m\lambda\rangle = t_1\langle h_{\beta_1}, \tau_1 (m+n)\lambda\rangle \in \bbZ.
    \]
    The remaining integrality constraints for $\pi'$ follow similarly.
    Assume $t_1 \ge m/(m+n)$, then we have $\pi''$:
    \begin{align*}
        t_1''\langle h_{\beta_1}, \tau_1 n\lambda\rangle
        &= \frac{(m+n)t_1 - m t_1'}{n}\langle h_{\beta_1}, \tau_1 n\lambda\rangle \\
        &= t_1\langle h_{\beta_1}, \tau_1 (m+n)\lambda\rangle - t_1'\langle h_{\beta_1}, \tau_1 m\lambda\rangle \in \bbZ.
    \end{align*}
    The remaining integrality constraints for $\pi''$ follow similarly.
    Thus, both $\pi'$ and $\pi''$ are LS paths of shape $m\lambda$ and $n\lambda$, respectively.

    Let $\psi\colon B(\pi_{(m+n)\lambda}) \to B(\pi_{m\lambda}) \otimes B(\pi_{n\lambda})$ denote the map $\pi \mapsto \pi' \otimes \pi''$.
    Since the embedding is unique and clearly $\pi_{(m+n)\lambda} \mapsto \pi_{m\lambda} \otimes \pi_{n\lambda}$, it is sufficient to show that $\psi$ commutes with $\te_i$ and $\tf_i$ for all $i \in I$.

    Fix $i\in I$.
    Let $t_{k^-}$ be the last time the global minimum of $\langle h_i, \pi(t)\rangle|_{[0,m/(m+n)]}$ is attained and let $t_{k^+}$ be the last time the global minimum of $\langle h_i, \pi(t)\rangle|_{[m/(m+n),1]}$ is attained.
    Now,
    \[
        \epsilon_i(\pi') = -\langle h_i, \pi(t_{k^-})\rangle
        \quad \text{and} \quad
        \varphi_i(\pi') = \langle h_i, \pi(m/(m+n)) - \pi(t_{k^-})\rangle
    \]
    and
    \[
        \epsilon_i(\pi'') = \langle h_i, \pi(m/(m+n)) - \pi(t_{k^+})\rangle
        \quad \text{and} \quad
        \varphi_i(\pi'') = \langle h_i, \pi(1) - \pi(t_{k^+})\rangle.
    \]
    From this we can clearly see that $\varphi_i(\pi') > \epsilon_i(\pi'')$ is equivalent to $\pi(t_{k^-}) < \pi(\tau_{k^+})$, and hence $\psi(\tf_i\pi) = \tf_i(\pi'\otimes \pi'')$ and $\psi(\te_i\pi) = \te_i(\pi'\otimes \pi'')$.
\end{proof}

From the explicit description of the embedding in Lemma \ref{lemma:crystal-cartan-embedding} we see that the partial order $\triangleright$ on $B(\pi_{(m+n)\lambda})$ is equivalent to the partial order $\triangleright$ on $B(\pi_{m\lambda} * \pi_{n\lambda})$.
Thus, we have the following result.

\begin{corollary}
    \label{corollary:left-action-straight-line-path-dilation}
    Let $\lambda$ be a dominant weight.
    Assume that Conjectures \ref{conjecture:f-action} and \ref{conjecture:e-action} hold for $V(\lambda)$ and the dominant path $\pi_\lambda$.
    Then, both conjectures also hold for $V(m\lambda)$ and the dominant path $\pi_{m\lambda}$ for all $m \ge 1$.
\end{corollary}

\begin{theorem}
    \label{theorem:base-change-path-monomial-basis-canonical-basis-classical-g2}
    Let $\frakg$ be a classical Lie algebra or of type $G_2$ and let $\lambda = m_1\varpi_1 + \dots + m_N\varpi_N$ be a dominant weight.
    Fix the dominant path $\eta := \pi_{m_1\varpi_1} * \dots * \pi_{m_N\varpi_N}$.
    Then, there is a unitriangular base change from path monomial basis $\calP(\eta)$ to the canonical basis of $V(\lambda)$ for any linearisation of $\trianglerighteq$ on $B(\eta)$.
\end{theorem}
\begin{proof}
    If $\lambda$ is a fundamental weight, Lemmas \ref{proposition:conjecture-classical-fundamental-weight} and \ref{proposition:conjecture-e-action-classical-fundamental-weight} imply that Conjectures \ref{conjecture:f-action} and \ref{conjecture:e-action} hold for $\pi_\lambda$.
    Applying Corollary \ref{corollary:left-action-straight-line-path-dilation} we see that both conjectures hold for $V(m\varpi_i)$ and the path $\pi_{m\varpi_i}$ for any $i\in I$ and $m \ge 0$.
    Using Proposition \ref{proposition:left-action-tensor-product}, we conclude that, in particular, Conjecture \ref{conjecture:f-action} holds for $V(\lambda)$ and $\eta$.
    Now, this theorem is a consequence of Proposition \ref{proposition:base-change-path-monomial-basis-canonical-basis}.
\end{proof}

\printbibliography

\end{document}